\documentclass[pdflatex,sn-nature]{sn-jnl}

\usepackage{graphicx}
\usepackage{amsmath,amssymb,amsfonts}
\usepackage{mathtools}
\usepackage{cases}
\usepackage[dvipsnames]{xcolor}
\usepackage{enumitem}
\usepackage{bm}
\usepackage{booktabs} 
\usepackage{tikz}
\usetikzlibrary{shapes.geometric, arrows, positioning}
\tikzstyle{every node}=[font=\small]
\tikzstyle{startstop} = [rectangle, rounded corners, minimum width=3cm,
    minimum height=1cm, text centered, draw=black, fill=gray!10]
\tikzstyle{process} = [rectangle, minimum width=3cm, minimum height=1cm,
    text centered, draw=black, fill=blue!5]
\tikzstyle{decision} = [diamond, aspect=2, minimum width=3cm, minimum height=1cm,
    text centered, draw=black, fill=yellow!10]
\tikzstyle{arrow} = [thick,->,>=stealth]

\begin{document}
% Replace the placeholders below with the final title, author details, abstract,
% and keywords before submission.
\title[Article title]{Let AEDs Move: Urban Mobility Enhanced Defibrillator Deployment
%with Evidence from New York City and Toronto
}

\author*[1]{\fnm{Bahar} \sur{D Viniche}}
\author[1]{\fnm{Sheng} \sur{Liu}}
\author[2]{\fnm{Nooshin} \sur{Salari}}

\affil*[1]{\orgname{Rotman School of Management, University of Toronto}}
\affil[1]{\orgname{Rotman School of Management, University of Toronto}}
\affil[2]{\orgname{DeGroote School of Business, McMaster University}}

\abstract{
Early defibrillation is critical to survival after out-of-hospital cardiac arrest (OHCA). However, public-access defibrillation programs rely largely on stationary automated external defibrillators (AEDs), whose fixed locations cannot adapt to spatial and temporal variation in OHCA demand and which a bystander must retrieve on foot.
We study a platform-enabled deployment model in which ride-hailing vehicles carry AEDs that can be dispatched directly to OHCA locations. Holding total AED supply fixed, we evaluate the full range of allocations between stationary and mobile deployment. Our analysis draws on public AED inventories, EMS-reported OHCA incidents, and ride-hailing mobility data from New York City and Toronto.
 %We reconstruct the spatiotemporal availability of the mobile fleet and evaluate AED response over the actual road network across the full range of \emph{hybrid ratios}, defined as the share of AED capacity allocated to the mobile layer.
% \textcolor{gray}{The results show that mobile AEDs can substantially improve both response times and reliability, but the benefits depend on how AED capacity is divided between stationary and mobile deployment and on which vehicles carry and respond with the devices.} 
Reallocating part of the supply to vehicles substantially shortens AED response times, but the benefit is non-monotonic in the mobile share. Under a specification common to both cities, mean improvement reaches approximately 2.9 minutes in each, yet the preferred allocation differs: it peaks with 55\% of AEDs mobile in Toronto, whereas New York City retains near-maximum gains across 40\%--80\%, with further variation across boroughs. The preferred allocation also depends on the response-time target, where coverage within 3 and 5 minutes favors greater mobility while 10-minute coverage favors retaining part of the stationary network.
Toronto's vehicle-level data further show that targeting AEDs to active vehicles is essential and with such targeting, restricting dispatch to idle vehicles, which disrupts no passenger trip, does not reduce the gain, and allowing any active equipped vehicle to respond raises it to 5.7 minutes under full mobilization.
Mobilizing AEDs can therefore materially improve access, and the stationary–mobile mix must be tailored to local coverage, fleet activity, and dispatch rules.
%and support city-specific deployment and operating policies rather than a universal stationary--mobile mix.
}


\maketitle

\section{Introduction and Background}
\label{sec:intro}

Out-of-hospital cardiac arrest (OHCA) is a significant public-health challenge, with more than 356{,}000 cases in the United States and approximately 60{,}000 in Canada each year \citep{ahaStats2022, heartstroke2024}.
Nearly 90\% of victims of OHCA die before reaching the hospital and less than 10\% survive to hospital discharge \citep{ahaStats2022, nichol2008regional, brooks2022optimizing}. 
%Most arrests occur at home or other non-public locations where bystanders are unlikely to have access to a defibrillator, and lay rescuers use AEDs in only 5.9--9\% of OHCAs \citep{weisfeldt2010, maltahansen2015}. 
Survival depends critically on the time to defibrillation, as treatment within a 3--5-minute window can achieve hospital-discharge survival rates above 50\%, while each additional minute without defibrillation reduces the likelihood of survival by about 10\%
\citep{larsen1993, valenzuela1997, valenzuela2000, merchant2020}.
Automated external defibrillators (AEDs) are one of the only effective methods for treating OHCAs, which enable untrained bystanders to deliver potentially lifesaving treatment before Emergency Medical Services (EMS) arrive. %When an OHCA incident happens, the nearby bystanders who observe the incident or are notified about it are instructed to retrieve the nearest AED and move to the OHCA site.
Therefore, cities have widely adopted public-access defibrillation (PAD) programs to make AEDs available in public spaces.
Yet %the effectiveness of these programs remains limited and
AEDs are retrieved and used in only a limited share of OHCAs, and just 2–4\% of patients receive defibrillation before professional responders arrive \citep{culley2004public, maltahansen2015}. 
% One central contributor is a spatiotemporal imbalance between AED supply and OHCA demand: arrests are broadly distributed across the city and across the day, whereas AEDs are fixed in place and their density is comparatively clustered.
This highlights a fundamental distinction between the availability of AEDs and timely access to defibrillation. Effective access depends jointly on the spatial and temporal proximity of an accessible device to the incident and on the time and resources required to bring it to the patient. Under conventional PAD practice, the latter involves a bystander traveling to retrieve an AED and returning with it to the arrest scene. %so that the time required to access the device reflects a round-trip journey rather than proximity alone.

Current PAD practice relies primarily on fixed-location AEDs placed in high-traffic public locations, including transportation hubs, airports, sports venues, and schools; we refer to these devices as \emph{stationary AEDs}.
A fundamental limitation of this design is a spatiotemporal mismatch between OHCA demand and AED supply.
Arrests are dispersed across the city, shift with population activity over the course of the day, with the majority occurring in residential settings \citep{weisfeldt2010}, whereas stationary AEDs are concentrated in public locations and are frequently housed within facilities whose access is restricted to operating hours.
For example, in Copenhagen, only 22.6\% of arrests occurred within 200 meters of a registered AED, and fewer than half of those devices were accessible at the time of the arrest \citep{karlsson2019}.
The resulting gap between registered and effective coverage is structural and reflects a broader challenge in urban infrastructure, in which public resources committed to fixed locations must serve demand that varies across space and time.
A body of research has focused on optimizing stationary AED placement to improve coverage of anticipated arrests and reduce the distance to the nearest device \citep{chan2016, chan2018}. 
However, even an optimized network remains spatially fixed, unable to adapt to changing patterns of emergency demand, and continues to rely on round-trip AED retrieval.
%The second direction introduces alternative mechanisms for bringing a device or a responder to the patient.
%Alternative response mechanisms have also been explored. 
Among alternative response mechanisms, drone-delivered AEDs fly the device directly to the patient and can shorten time to defibrillation \citep{pulver2016, boutilier2017, claesson2017}. Drones, however, are dedicated assets that operate from fixed bases and require purpose-built infrastructure and centralized dispatch. Community-responder programs instead use smartphone alerts to direct nearby volunteers to the patient \citep{andelius2020, pulsepoint}, expanding the pool of potential responders without changing the underlying retrieval process or AED network. %\textcolor{gray}{These limitations motivate an alternative approach that combines dynamic AED network with direct delivery, without requiring a separate emergency-delivery fleet.}

We propose a platform-enabled AED deployment model that introduces a mobile layer to stationary coverage by equipping a subset of ride-hailing vehicles with AEDs and integrating them into emergency dispatch. Ride-hailing vehicles circulate across cities as part of their normal operations, with their spatial distribution shifting throughout the day in response to commuting flows and other patterns of urban activity \citep{gong2015, mckenzie2020, rahimi2021}. Rather than remaining at predetermined locations, the AEDs they carry move with flows and potentially bring coverage into closer alignment with the spatial and temporal distribution of emergency demand.
When an OHCA occurs, an equipped vehicle with mobile AED can be dispatched directly to the patient, replacing a bystander's round-trip retrieval on foot with one-way delivery by vehicle. %\textcolor{gray}{This model thus leverages existing fleet movements to provide dynamic AED coverage and direct delivery.}

Practical experience provides initial evidence of the potential of mobile AED deployment. In Singapore, the AED-on-Wheels pilot equipped 120 taxis with AEDs and alerted drivers to nearby OHCAs. Between 2015 and 2017, participating taxis accepted 374 emergency alerts and arrived on scene in an average of 6 minutes and 22 seconds, faster than the reported ambulance response time \citep{ng2021}. 
Policy interest is also growing in the United Kingdom, where the JumpStart Campaign advocates equipping all newly sold cars with AEDs, so that privately owned vehicles can form a mobile network of defibrillators \citep{jumpstart2026}. %\textcolor{gray}{often parked near the homes where most arrests occur,}
In parallel, the driving theory test, taken about 2.4 million times a year, now includes questions on AED use, so that drivers are prepared to act as first responders \citep{dvsa2025theory}.
While these initiatives demonstrate operational feasibility and policy interest, they do not provide a systematic evaluation of mobile AED performance or guidance on how mobile deployment should be integrated with existing stationary coverage.
Related research has examined AEDs carried by public buses in Pittsburgh, showing their potential to extend coverage alongside stationary AEDs \citep{hajari2020simulating}. However, buses remained on predetermined routes and schedules, unable to reposition dynamically in response to changes in urban activity, and AED access still required bystander retrieval.
In concurrent work, \citep{lin2026} developed a stochastic analytical model comparing stationary and mobile AED deployment on a stylized circular road, characterizing the role of population, budget, and fleet idleness. Such abstractions provide valuable analytical insights, but the stylized setting abstracts from the road-network structure and actual vehicle movements.
These features are central to our evaluation of mobile AED systems since their effectiveness critically depends on where vehicles are located when an OHCA occurs and how quickly they can reach the patient.
%for final verison: Related analytical work examines stationary and mobile AED deployment using a stochastic model on a stylized circular road, providing complementary theoretical insights \citep{lin2026}.

% We refer to the deployment concept as \emph{AED-on-Wheels}, and to AEDs deployed through the mobility fleet as \emph{mobile AEDs}.
Do mobile AEDs improve timely access in real urban settings, and by how much? And how should a city's AED supply be divided between stationary and mobile deployment?
We address these questions by holding total AED supply fixed and comparing alternative allocations of the same devices, which isolates the value of mobility from that of adding devices.
The full stationary network serves as the status-quo benchmark, from which we progressively reallocate a share of the AED supply to ride-hailing vehicles, ranging from a fully stationary to a fully mobile network. Intermediate allocations form \emph{hybrid policies}, in which stationary and mobile AEDs operate jointly, and we refer to the share of supply allocated to vehicles as the \emph{hybrid ratio}.
This reallocation entails a trade-off: each additional equipped vehicle extends mobile coverage, while each device withdrawn from a stationary site may erode existing coverage. The preferred allocation therefore depends on how effectively the mobile fleet compensates for the stationary coverage it replaces.
Using EMS-reported OHCA incidents, public AED inventories, and ride-hailing mobility data from Toronto and New York City, we provide, to our knowledge, the first empirical evaluation of reallocating a fixed AED supply between stationary and mobile deployment.

We evaluate each allocation by its mean response-time improvement over the stationary benchmark, the distribution of gains and losses across OHCA events, and the reliability of timely access, computing response times over the road network.
%We  under assumptions that deliberately favor the stationary benchmark.
Our cross-city analysis shows that the benefit of mobilization can be non-monotonic, with an intermediate stationary--mobile allocation outperforming both fully stationary and fully mobile deployment.
The location and breadth of this advantage vary between and within cities, reflecting the interplay of road-network structure, existing AED coverage, the geography of cardiac arrests, and fleet mobility.
Toronto's vehicle-level data, which distinguish idle, en-route, and occupied phases, further show that whether and where this optimum arises depends on which vehicles carry AEDs and which may respond.
The appropriate stationary--mobile mix is therefore not universal but depends on local coverage, fleet activity, and dispatch rules.
%: it emerges when the dispatchable equipped fleet is too small to replace the withdrawn stationary coverage, whereas under the most permissive dispatch rule the improvement rises to 5.7 minutes under full mobilization.
%\textcolor{gray}{These results quantify the expected reduction in OHCA response times and the associated public-health implications of mobile AED deployment relative to the current stationary AED infrastructure.}

% Assumptions (for authors): baseline treats every registered AED as accessible at all hours and
% omits recognition/retrieval delay (favors stationary); mobile leg omits notify-accept-reroute
% latency (favors mobile). Net direction is not established, so we do not claim the setup favors stationary.

\section{Problem Setting}

AED-on-Wheels deploys mobile AEDs through an existing urban mobility platform, such as a ride-hailing fleet. Under the conventional stationary system, AEDs remain at fixed locations and, when a second bystander is available, the bystander must travel to the nearest stationary AED and return with it to the OHCA patient. Under AED-on-Wheels, mobile AEDs move throughout the city as part of normal fleet operations and can be dispatched directly to an OHCA location when needed.

The potential advantage of mobile AEDs arises from two features. First, stationary AED retrieval requires a two-way trip by a bystander, whereas a mobile AED travels directly toward the patient. 
Second, the case for using ride-hailing vehicles as mobile AED carriers is also based on the conjecture that the distribution of ride-hailing vehicles better aligns with the OHCA distribution than that of stationary AEDs. 
Ride-hailing fleets continuously reposition throughout the city as travel demand and urban activity change. The spatial distribution of mobile AEDs may therefore align more closely with where people are and, consequently, where potential emergency events are likely to occur.
Moreover, vehicle locations adapt naturally to temporal variations such as rush hours, weekends, and night periods, maintaining an approximate proportionality between vehicle density and OHCA likelihood.
This spatiotemporal alignment helps reduce response times. We formalize this intuition analytically below to isolate this mechanism.

\subsection{Why AED-on-Wheels May Help}
Consider a service region denoted by $\Omega\subset\mathbb{R}^2$ with area $A=|\Omega|$.  
Let $g(x)$ denote the spatial intensity of OHCA events and $G \;=\; \int_{\Omega} g(x)\,dx$ denote the total event rate. 
We assume the system has a fixed total of $N$ AEDs and compare two pure deployment rules using the same capacity: stationary AED deployment and fully mobile AED deployment.

The following spatial models are analytical benchmarks rather than descriptions of the observed stationary and mobile AED distributions in our empirical settings. They isolate the potential value of aligning AED supply with the spatial distribution of emergency demand.

%\subsubsection{Mean Distance- Response Time}
\noindent\textbf{Scenario 1 (uniform AED distribution)}  
Assume AEDs are placed uniformly at fixed locations (stationary); then the AED intensity 
$\lambda_K^{(1)} = \frac{N}{A}$ is constant. In this case, the expected nearest-AED distance is
\[
\mathbb{E}[R\mid x] \;=\; \frac{1}{2\sqrt{\lambda_K^{(1)}}},
\]
thus, the (event-weighted) mean nearest distance under Scenario\,1 is
\[
\mathbb{E}[R]_{1}
\;=\;
\frac{1}{2\sqrt{\lambda_K^{(1)}}}
\;=\;
\frac{1}{2}\sqrt{\frac{A}{N}}.
\]

\noindent\textbf{Scenario 2 (AEDs distributed proportional to events)}  
If we can place AEDs according to the event distribution, i.e., the AED density satisfies
\[
\lambda_K^{(2)}(x)=\alpha\,g(x),\qquad \alpha=\frac{N}{G}, 
\]
so that the total number of AEDs integrates to $N$.  Using the local approximation
$\mathbb{E}[R\mid x]\approx 1/(2\sqrt{\lambda_K^{(2)}(x)})$ and averaging over event locations
with weight $g(x)/G$ yields
\[
\mathbb{E}[R]_{2}\;\approx\;
\int_{\Omega}\frac{1}{2\sqrt{\alpha\,g(x)}}\cdot\frac{g(x)}{G}\,dx
\;=\;
\frac{1}{2\sqrt{\alpha}}\cdot\frac{1}{G}\int_{\Omega}\sqrt{g(x)}\,dx.
\]
Substituting $\alpha=N/G$, the (event-weighted) mean nearest distance under Scenario\,2 is
\[
\mathbb{E}[R]_{2}
\;=\;
\frac{1}{2}\,\frac{\displaystyle\int_{\Omega}\sqrt{g(x)}\,dx}{\sqrt{G\,N}}.
\]

\noindent\textbf{Comparison:}  The ratio of the response distances under the two scenarios is
\[
\frac{\mathbb{E}[R]_{2}}{\mathbb{E}[R]_{1}}
\;=\;
\frac{\dfrac{1}{2}\,\dfrac{\int_{\Omega}\sqrt{g(x)}\,dx}{\sqrt{G\,N}}}
{\dfrac{1}{2}\sqrt{\dfrac{A}{N}}}
\;=\;
\frac{\displaystyle\int_{\Omega}\sqrt{g(x)}\,dx}{\sqrt{G\,A}}.
\]
By the Cauchy--Schwarz inequality,
$\int_{\Omega}\sqrt{g(x)}\,dx \le \sqrt{\int_{\Omega}1^2 dx}\,\sqrt{\int_{\Omega} g(x)\,dx}
= \sqrt{A\,G},$
so we have
\[
\mathbb{E}[R]_{2} \le \mathbb{E}[R]_{1},
\]
with equality if $g(x)$ is constant (i.e., events are uniform). The analysis implies that if AEDs are allocated proportionally to event intensity, the event-weighted mean nearest distance cannot exceed that of uniform placement; proportional allocation strictly reduces the mean distance whenever $g$ is nonuniform. Furthermore, the improvement magnitude is governed by the factor $\int_{\Omega}\sqrt{g(x)}\,dx/\sqrt{G\,A}$: the more spatially concentrated (``peaked'') $g$ is, the smaller this factor and the larger the relative gain from proportional placement. %(iii) In the uniform-$g$ limit, the two strategies coincide.
%These closed-form relations provide a straightforward metric to assess the potential benefits of mobile AEDs in response to observed spatial heterogeneity in OHCA incidence.

% \textcolor{Green}{The analysis is also extendable to the hybrid setting, where a portion of AEDs (stationary AEDs) are uniformly distributed, and a portion are mobile AEDs with the correlated distribution.}

We can also analyze reliability performance when the policymaker aims to meet a response-time target with high probability. Given $\lambda_K^{(2)}(x)=\alpha\,g(x)$, suppose that the number of mobile AEDs in a neighborhood follows a spatial Poisson process with density $\lambda_K^{(2)}(x)=\alpha\,g(x)$, the probability that at least one AED is available in a disk of radius $r_m$ is
\begin{equation*}
R_{\mathrm{mob}}(r_m)
= 1 - \mathbb{E}\!\left[\exp\!\Big(-\pi r_m^2\,\alpha\,g(X)\Big)\right],
\qquad X\sim g(x)/G,
\end{equation*}
where the expectation is taken over incident locations $X$ with size-biased density $g(x)/G$ (incidents occur more often where $g$ is larger). Applying Jensen's inequality ($e^{-x}$ is convex) gives
\begin{equation*}
R_{\mathrm{mob}}(r_m) \le 1 - \exp\!\Big(-\pi r_m^2\,\alpha\,\overline g_*\Big),
\quad
\overline g_*=\frac{1}{G}\int_\Omega g(x)^2\,dx.
\end{equation*}
Thus, the availability probability improves exponentially in $\alpha$ and in the size-biased mean $\overline{g}_*$, which increases with spatial concentration of OHCA risk. Building on this analytical intuition, we perform large-scale numerical analysis using real-world mobility data from New York City and Toronto to assess the performance of Mobile-AED policies in different forms.

\section{AED-on-Wheels Policy Evaluation}

In this section, we evaluate the proposed AED-on-Wheels policy using real-world data from Toronto and New York City. For each city, we assemble three layers of data: (i) the OHCA-proxy incident set derived from EMS records, (ii) the public inventory of stationary AEDs, and (iii) the mobility layer of vehicle trips used to represent the spatiotemporal availability of potential mobile AED carriers.
% We report the data construction and processing choices required to interpret the policy analysis.
% while construction, processing procedures, geocoding, routing, and additional descriptive analyses of each layer are detailed in Appendix. 

% \subsection{Policy Evaluation}
\label{sec:policy-eval}

We compare three deployment policies under a fixed AED supply. The \emph{stationary} policy is the status quo, in which AEDs remain at fixed locations and a bystander travels to the nearest AED and brings it back to the OHCA scene. The \emph{mobile} policy places AEDs on circulating ride-hailing vehicles, and the nearest available mobile AED is dispatched directly to the OHCA. The \emph{hybrid} policy reallocates a share of the existing stationary AED inventory to the mobile fleet while retaining the remaining capacity at fixed locations. For each OHCA event, the response time is determined by the faster of a bystander's two-way retrieval of a remaining stationary AED and a vehicle's one-way delivery. We index these policies by the \emph{hybrid ratio}, the share of total AED capacity allocated to vehicles, ranging from $0\%$ (fully stationary) to $100\%$ (fully mobile).
% I) The status quo Stationary AED deployment policy where AEDs are solely deployed in fixed locations which provides a baseline for understanding the existing accessibility and coverage of AEDs,
% II) Mobile AED policy where AEDs are only mounted on taxis and ride-hailing vehicles, offering the potential for a more responsive and timely AED deployment, and
% III) Hybrid Policy where existing stationary AEDs continue to serve as a foundation, while selected taxis and ride-hailing vehicles are equipped with mobile AEDs to provide comprehensive coverage and enhancing both the spatial reach and temporal responsiveness of the system.
For each hybrid ratio, we construct a counterfactual stationary network by withdrawing the corresponding number of devices from the observed AED inventory. We evaluate ten independent mobility samples at each hybrid ratio to account for random vehicle selection. Performance is measured as the improvement relative to the full stationary network, aggregated across OHCA events and averaged across the mobility samples. In computing response times, we assume that a second bystander is always available, that stationary retrieval incurs no delay for locating or unlocking the AED, and that vehicles travel at a speed representative of heavy traffic. The construction is detailed in Methods.

% \textcolor{blue}{A ride-hailing vehicle that is logged in to the platform is, at any instant, in one of three \emph{activity phases} (Figure~\ref{fig:toronto_phases}): \emph{idle} (P1), available on the platform without a passenger or an assigned request; \emph{en route} (P2), travelling to pick up an assigned passenger; and \emph{occupied} (P3), carrying a passenger between pickup and drop-off. An equipped vehicle carries its AED through all three phases, so the phases do not change which vehicles hold AEDs, only which of them may respond at the instant of an OHCA. We refer to the set of phases in which an equipped vehicle may be dispatched as the \emph{dispatch rule}. Under the \emph{occupied} rule only vehicles in P3 respond, under the \emph{idle} rule only vehicles in P1, and under the \emph{all-phase} rule any active equipped vehicle responds.}

\subsection{Data}
\label{sec:data}

\subsubsection{Toronto}
\label{sec:toronto}

For the City of Toronto, we analyze AED deployment policies over the full calendar year 2024 using OHCA-proxy from the Toronto Paramedic Services Incident Data \citep{TorontoParamedic2025}, the City of Toronto's public AED inventory \citep{TorontoAED2025}, and twelve months of Vehicle-for-Hire (VFH) records from the two major ride-hailing platforms operating in the city.
%, obtained directly from the City of Toronto.

% The OHCA-proxy events are medical incidents dispatched at the highest priority level in the EMS records, which is assigned to suspected cardiac arrest.
% These incidents are then downsampled to the annual cardiac-arrest count implied by Toronto Paramedic Services reports, yielding $4{,}817$ events.
The The OHCA-proxy events are medical incidents dispatched at the highest priority level, which is reserved for suspected cardiac and respiratory arrest and as a proxy for OHCA; these incidents are downsampled to the annual cardiac-arrest count reported by Toronto Paramedic Services ($4{,}817$ events).

The stationary layer is the City of Toronto public AED inventory \citep{TorontoAED2025}. After geocoding, and excluding devices located outside the city boundaries, the inventory contains $1291$ AEDs at $902$ distinct locations. Of these $1291$ devices, $389$ ($30.1\%$) are additional devices at locations that already host one.
%, compared with $33.8\%$ in New York City.

\begin{figure}[!t]
    \centering
    \includegraphics[width=0.8\linewidth]{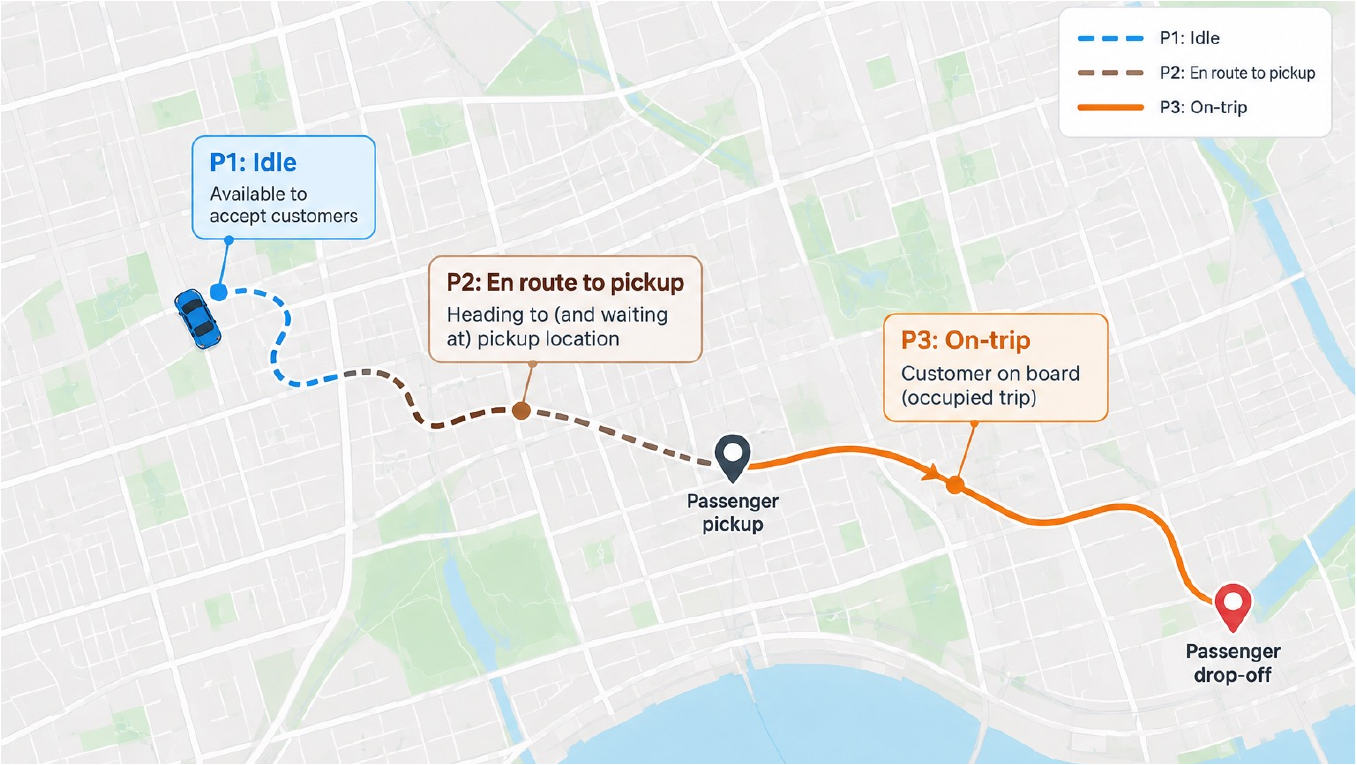}
    \caption{Schematic illustration of the three ride-hailing activity phases: idle (P1), en route to pickup (P2), and occupied trip (P3).}
    \label{fig:toronto_phases}
\end{figure}

The mobility layer consists of twelve months of VFH records, which are not publicly available. 
These data follow each vehicle through three \emph{activity phases}: \emph{idle} (P1), available on the platform without a passenger or an assigned request; \emph{en route} (P2), travelling to pick up an assigned passenger; and \emph{occupied} (P3), carrying a passenger between pickup and drop-off (Figure~\ref{fig:toronto_phases}).

% \textcolor{gray}{The three datasets are released at different spatial resolutions, so we make simple location assumptions where exact coordinates are not available; the construction of each layer is described in Appendix~\ref{app:data}.}
% [MOVED TO METHODS] Incidents are reported by postal area (Forward Sortation Area, FSA), and each OHCA event is placed at the population-weighted centroid of its FSA, the Toronto analogue of the ZIP-code centroids used for New York City. AEDs are placed at their geocoded addresses, and vehicles are positioned on the road network from the location data provided, so that each vehicle can be located at any instant. Further detail on the vehicle location data is withheld for privacy reasons.
% [MOVED TO METHODS] footnote on the equipping pool: "Drawn from the vehicles with at least five occupied trips in that month, with probability proportional to their activity."

Two features of the Toronto VFH data enable a detailed representation of mobile AED deployment and offer additional insight into how the policy would operate. First, the VFH records carry a persistent vehicle identifier, so mobile AEDs are assigned to vehicles rather than to individual trips. At each hybrid ratio, the equipped set is drawn from the vehicles active on the platform in that month, with each vehicle’s chance of selection proportional to its activity on the platform (\emph{activity-weighted equipping}). %The set’s size equals the reallocated share of the $1{,}291$ AEDs, nested across ratios, and drawn independently in each of ten mobility samples. 
Second, because the records distinguish the activity phase of each vehicle, we can vary which equipped vehicles are allowed to respond to an OHCA event. This separates two questions of \textit{whether ride-hailing vehicles are close enough to an OHCA to improve response}, and \textit{whether the vehicles that are close could plausibly be diverted to it}. A vehicle carrying a passenger could be a less plausible responder, while an idle vehicle can be diverted without disrupting a trip. These are detailed in Section \ref{sec:toronto-availability}.

\subsubsection{New York City}
\label{sec:nyc}
For New York City, we analyze AED deployment policies over the full calendar year 2024 using three datasets: OHCA-proxy events from the NYC EMS Incident Dispatch data \citep{FDNYEMSDispatch}, the NYC Open Data public AED inventory \citep{NYCAED2026}, and the High-Volume For-Hire Vehicle (FHV) trip records published by the NYC Taxi \& Limousine Commission (TLC) \cite{NYCTLCTripData}, which cover every trip served by the major ride-hailing platforms in the city, including Uber and Lyft. %are publicly available,
The OHCA-proxy events are EMS incidents dispatched and closed as cardiac arrest, which yields $10{,}082$ events over 2024. 
The AED inventory contains $11{,}138$ AEDs at $7{,}373$ distinct locations, of which $33.8\%$ are additional devices at locations that already host one.

In contrast to Toronto's vehicle-level records, the NYC mobility data are less granular in three respects. First, the records describe trips rather than vehicles and carry no vehicle identifier, so mobile AEDs are assigned to a sample of trips at each hybrid ratio rather than to individual vehicles. Second, only the occupied phase between pickup and drop-off is observed, so an equipped vehicle can respond only while carrying a passenger. Third, trip origins and destinations are reported at the level of taxi zones. 
The New York City analysis therefore corresponds to the occupied dispatch rule of Section~\ref{sec:toronto-hybrid}. 
Because AEDs are assigned to trips, the mobile layer follows occupied ride-hailing activity in proportion to its volume, which is closest in spirit to the activity-weighted equipping used in Toronto.
The New York City data compensate in scale.
The FHV records cover between $15$ and $20$ million trips per month across all five boroughs. This extensive geographic coverage allows us to go beyond citywide performance and examine how mobile AED deployment and the preferred stationary--mobile allocation performs across distinct urban settings within a single city (Appendix~\ref{app:nyc-borough}).

\subsection{Stationary Coverage}
\label{sec:baseline}

\begin{figure}[!b]
\centering
\includegraphics[width=\textwidth]{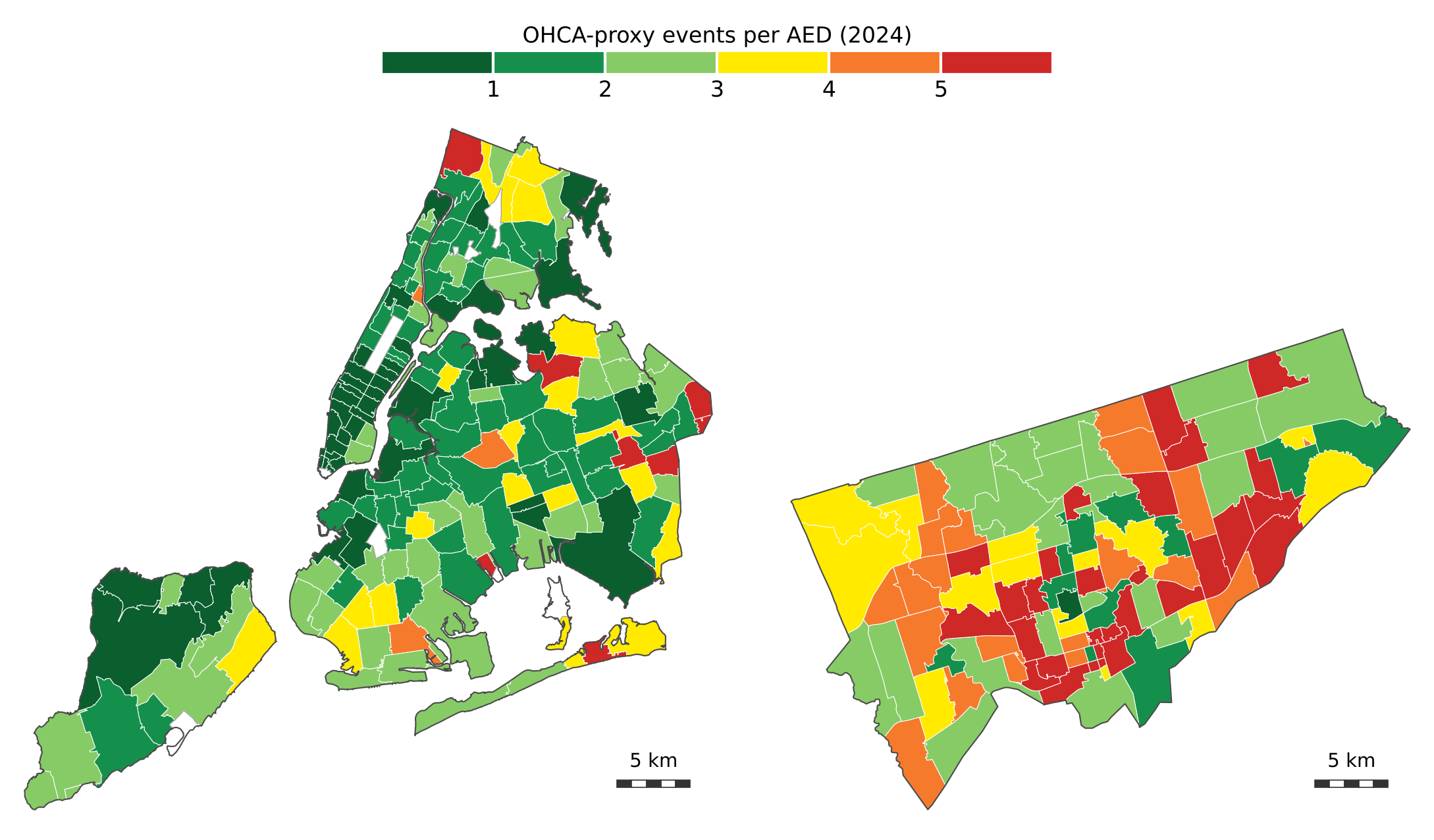}
\par\smallskip
\begin{minipage}[t]{0.524\textwidth}\centering{\small (a) Spatial distribution in NYC}\end{minipage}%
\begin{minipage}[t]{0.476\textwidth}\centering{\small (b) Spatial distribution in Toronto}\end{minipage}
\par\bigskip
\begin{minipage}[t]{0.49\textwidth}
    \centering
    \includegraphics[width=\textwidth]{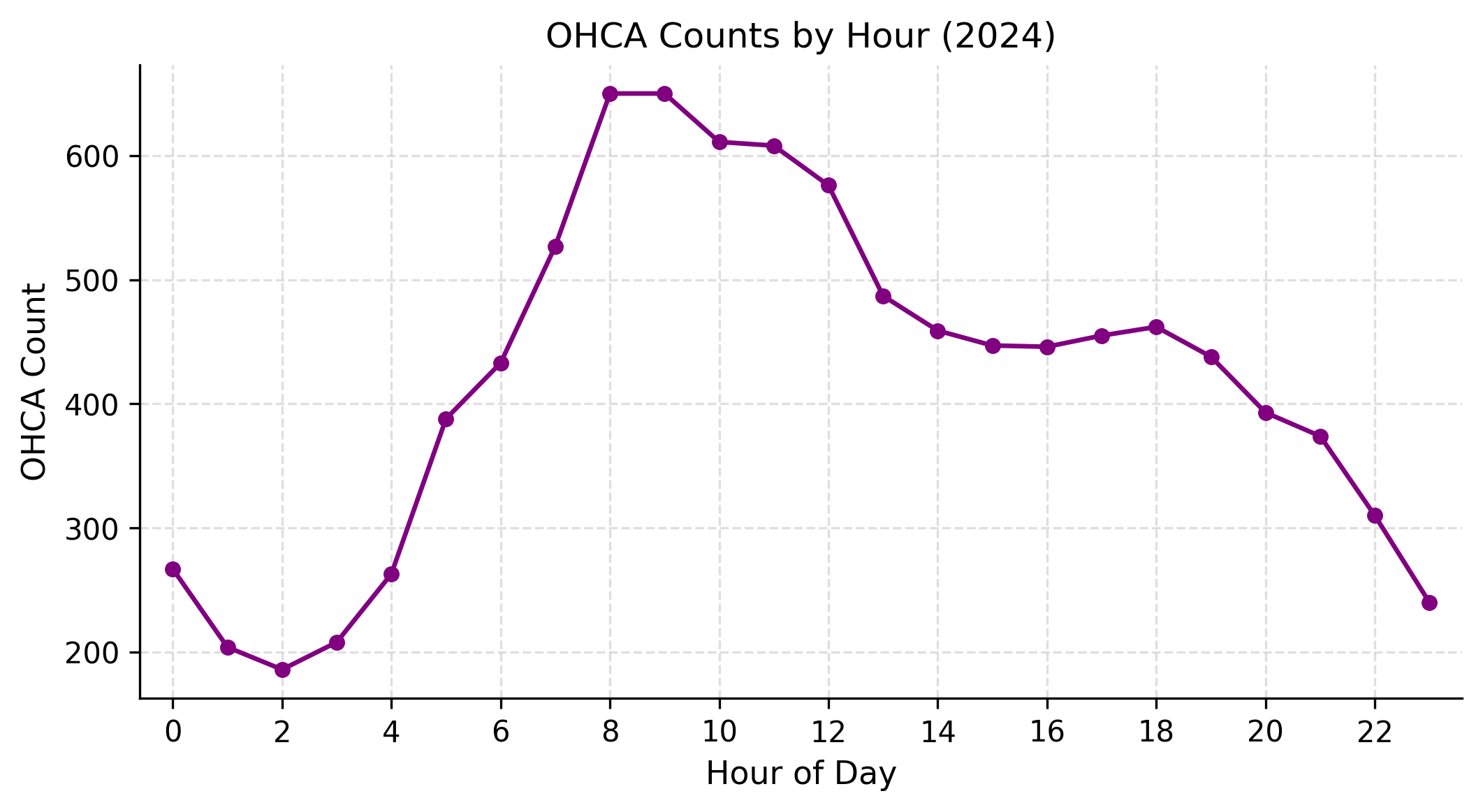}
    \par\smallskip{\small (c) Temporal distribution of OHCA events - NYC}
\end{minipage}\hfill
\begin{minipage}[t]{0.49\textwidth}
    \centering
    \includegraphics[width=\textwidth]{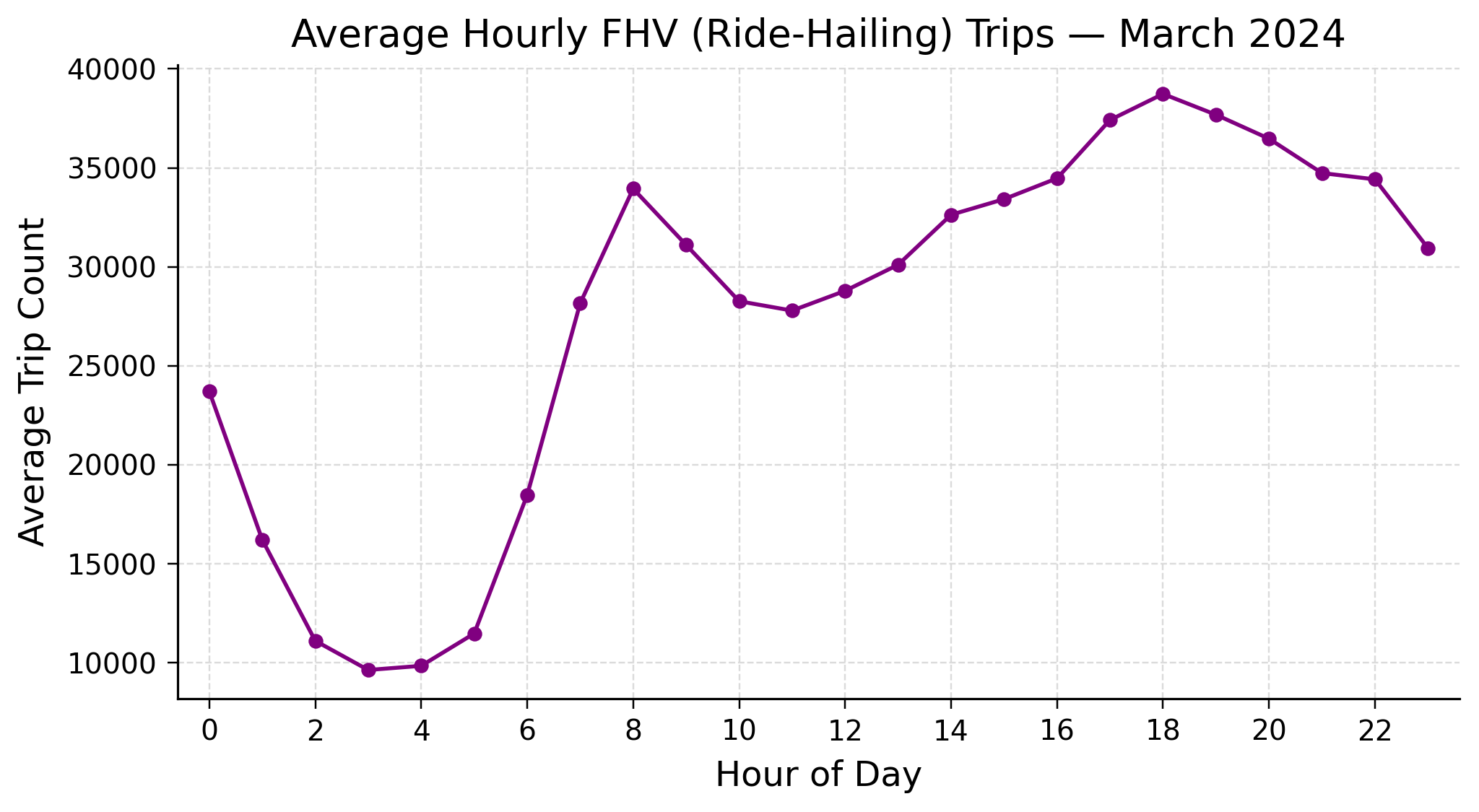}
    \par\smallskip{\small (d) Temporal distribution of ride-hailing (for-hire-vehicle) trips - NYC}
\end{minipage}
\caption{Spatial and temporal patterns relevant to AED accessibility. (a, b) The OHCA-to-AED ratio by spatial unit (ZIP code in New York City, FSA in Toronto), on a common scale, shows substantial spatial heterogeneity in AED-to-OHCA ratio in both cities. (c) OHCA-proxy events in New York City are distributed unevenly across the hours of the day. (d) Average for-hire-vehicle trip counts in New York City in March 2024, a representative month, throughout the hours of the day.}
\label{fig:nyc_align}
\end{figure}

We first characterize the existing stationary network and how well it aligns with the demand it is meant to serve. Figure~\ref{fig:nyc_align} a,b spatial distribution of OHCA-proxy events relative to stationary AED supply in both cities. In New York City, the OHCA-to-AED ratio rises outside the dense core, whereas in Toronto it remains high across much of the city, including downtown. The temporal patterns in New York City further highlights the variation in OHCA demand and ride-hailing activity throughout the day (Figure~\ref{fig:nyc_align} c,d). %These patterns point to the spatial and temporal misalignment that motivates this study, between where and when AEDs are fixed and where and when OHCA events occur, which a mobility-enhanced deployment could reduce in both dimensions.
 
% Toronto's stationary network is markedly sparser relative to OHCA demand than New York City's. Under the full stationary network (hybrid ratio $0\%$), the mean on-foot network distance from an OHCA event to its nearest AED is $503$~m and the median is $457$~m, compared with $316$~m and $271$~m in New York City over the same year (Table~\ref{tab:borough_baseline}). The corresponding two-way retrieval time averages $12.9$~minutes. Only $2.2\%$, $4.7\%$, and $41.8\%$ of OHCA events are reachable within $3$, $5$, and $10$~minutes, against $10.6\%$, $32.3\%$, and $68.7\%$ in New York City. 

% Under the full stationary network, the mean on-foot network distance from an OHCA event to its nearest AED is $316$~m and the median is $271$~m, and the two-way retrieval takes $8.1$~minutes on average; $10.6\%$, $32.3\%$, and $68.7\%$ of OHCA events are reachable within $3$, $5$, and $10$~minutes. 
% Within the city, this coverage is itself not uniform across boroughs (Table~\ref{tab:borough_baseline}): Manhattan reaches $93.2\%$ of its events within $10$~minutes, whereas Richmond/Staten Island reaches $36.4\%$.
% These differences establish heterogeneous starting conditions across the city, where the value of reallocating AED capacity to the mobile network depends on existing stationary coverage and the local availability of the mobility fleet.

Toronto’s stationary network is markedly sparser relative to OHCA demand than New York City's. Under the full stationary network, the mean on-foot network distance from an OHCA event to its nearest AED is 503 m in Toronto, compared with 316 m in New York City, with corresponding medians of 457 m and 271 m. The average two-way retrieval time is 12.9 minutes in Toronto and 8.1 minutes in New York City. The difference is particularly pronounced at the critical 5-minute threshold: only 4.7\% of OHCA events in Toronto are reachable, compared with 32.3\% in New York City. Within 10 minutes, the corresponding shares are 41.8\% and 68.7\%.
Within New York City, coverage is heterogeneous across boroughs (Table~\ref{tab:borough_baseline}): Manhattan reaches 93.2\% of its events within 10 minutes, whereas Richmond/Staten Island reaches 36.4\%. These differences establish distinct starting conditions, where the value of reallocating AED capacity to the mobile network depends on existing stationary coverage and the local availability of the mobility fleet.

\subsection{Hybrid Mobile-AED Deployment Policy }
\label{sec:toronto-hybrid}

The New York City trip records cover the occupied phase alone. To place the two cities on the same basis, we allow an equipped vehicle to respond while carrying a passenger (P3), which we call the \emph{occupied dispatch rule}. We then extend the Toronto analysis to other dispatch rules and compare the policy implications in Section~\ref{sec:toronto-availability}. %That extension shows that restricting dispatch to idle vehicles, which avoids interrupting passenger trips, yields similar or larger mean improvement, so the occupied rule does not overstate what such a policy can achieve. 
We evaluate the hybrid policy on the $4{,}817$ OHCA-proxy events of 2024 in Toronto and the $10{,}082$ in New York City across the full range of hybrid ratios, following the construction of Section~\ref{sec:policy-eval}, with outcomes pooled over the twelve months and averaged across the ten mobility samples.

\subsubsection{Mean Response-Time Improvement}
\label{sec:hybrid}

\begin{figure}[!b]
\centering
\begin{minipage}[t]{0.49\textwidth}
    \centering
    \includegraphics[width=\textwidth]{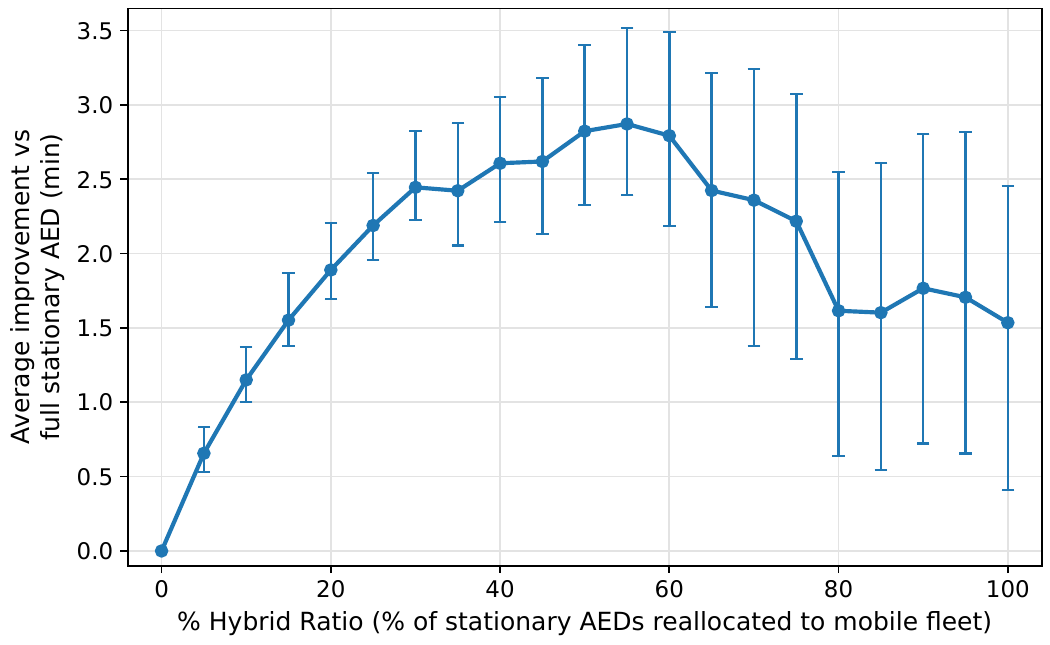}
    \par\smallskip{\small (a) Toronto, 2024}
\end{minipage}\hfill
\begin{minipage}[t]{0.49\textwidth}
    \centering
    \includegraphics[width=\textwidth]{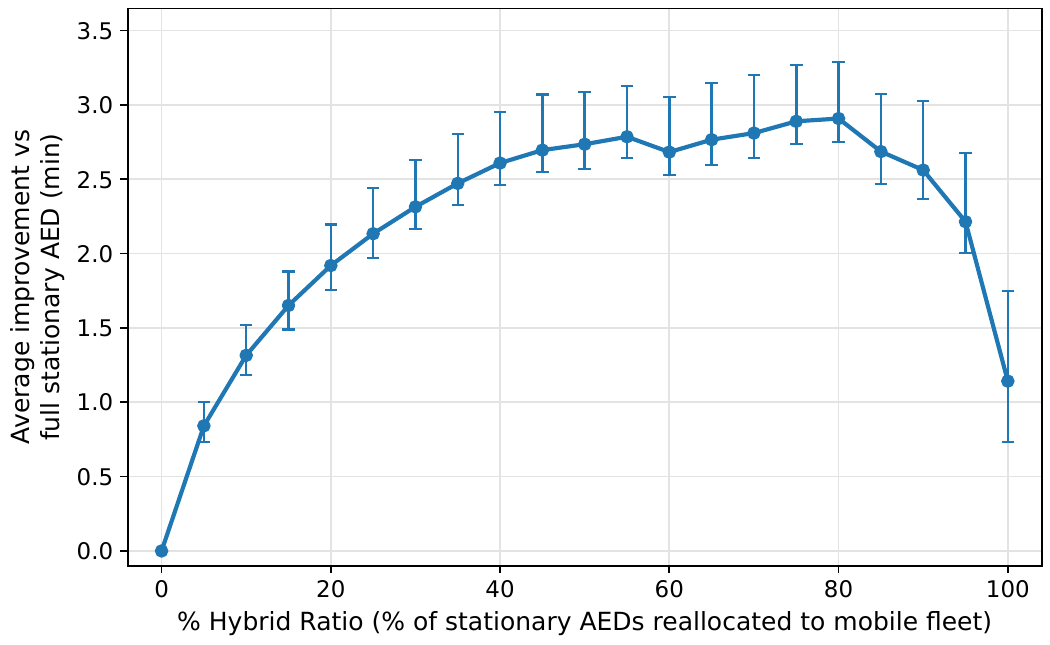}
    \par\smallskip{\small (b) New York City, 2024}
\end{minipage}
\caption{Mean response-time improvement relative to the full stationary network, as a function of the hybrid ratio, under the occupied dispatch rule (2024): (a) Toronto, (b) New York City. Points are averages over the year and the ten mobility samples; bars show the range across months.}
\label{fig:ratio_improvement}
\end{figure}

Moving part of the AED supply from fixed locations onto ride-hailing vehicles shortens the average time for an AED to reach an OHCA at every hybrid ratio. Figure~\ref{fig:ratio_improvement} shows that in both cities the benefit follows the same pattern as the hybrid ratio increases: mean improvement rises steeply as the first AEDs are moved to vehicles, levels off at intermediate hybrid ratios, and declines as the network approaches a fully mobile one. The pattern reflects a trade-off inherent in reallocating a fixed AED supply, since every device placed on a vehicle extends mobile coverage while withdrawing stationary coverage. When few AEDs are mobile, each additional equipped vehicle reaches events that no nearby fixed device serves, and the added mobile coverage dominates. As reallocation continues, the loss of stationary coverage accelerates, and at high hybrid ratios this exceeds what the mobile fleet can replace. As a result, an intermediate allocation outperforms both the full stationary network and a fully mobile one, although even a fully mobile network still improves on the stationary benchmark. The steep initial rise also reflects the reallocation rule. Since existing AED networks may concentrate multiple devices at few locations, %reference to Tim's paper
we prioritize reallocating these additional devices, so no AED \textit{location} is actually withdrawn at low hybrid ratios.

\begin{figure}[htbp]
\centering
\begin{minipage}[t]{0.49\textwidth}
    \centering
    \includegraphics[width=\textwidth]{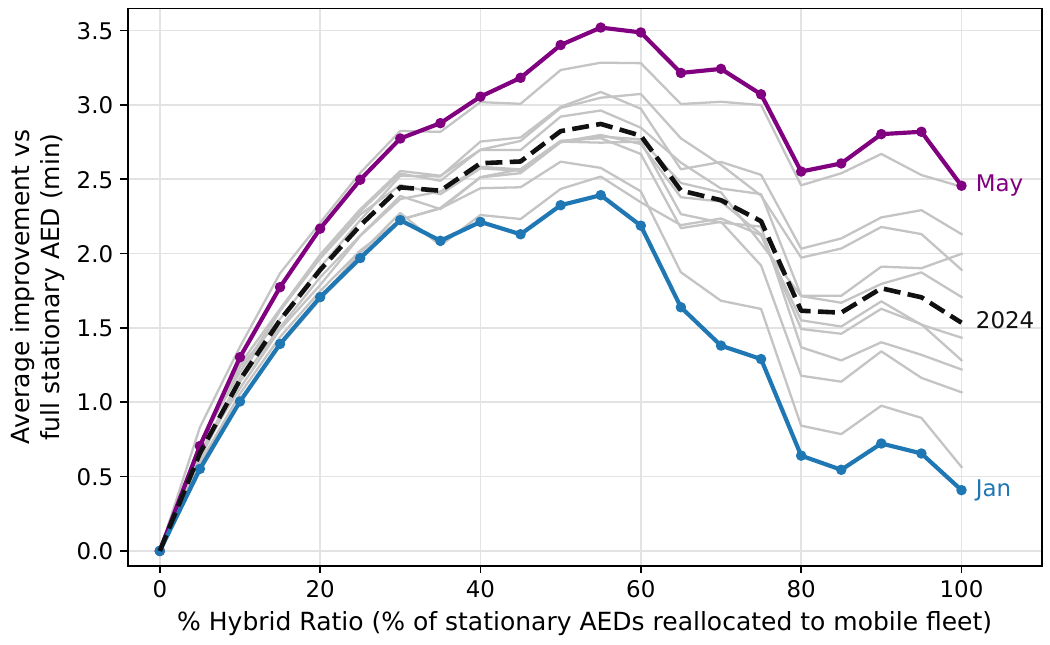}
    \par{\small (a) Mean improvement by month}
\end{minipage}\hfill
\begin{minipage}[t]{0.49\textwidth}
    \centering
    \includegraphics[width=\textwidth]{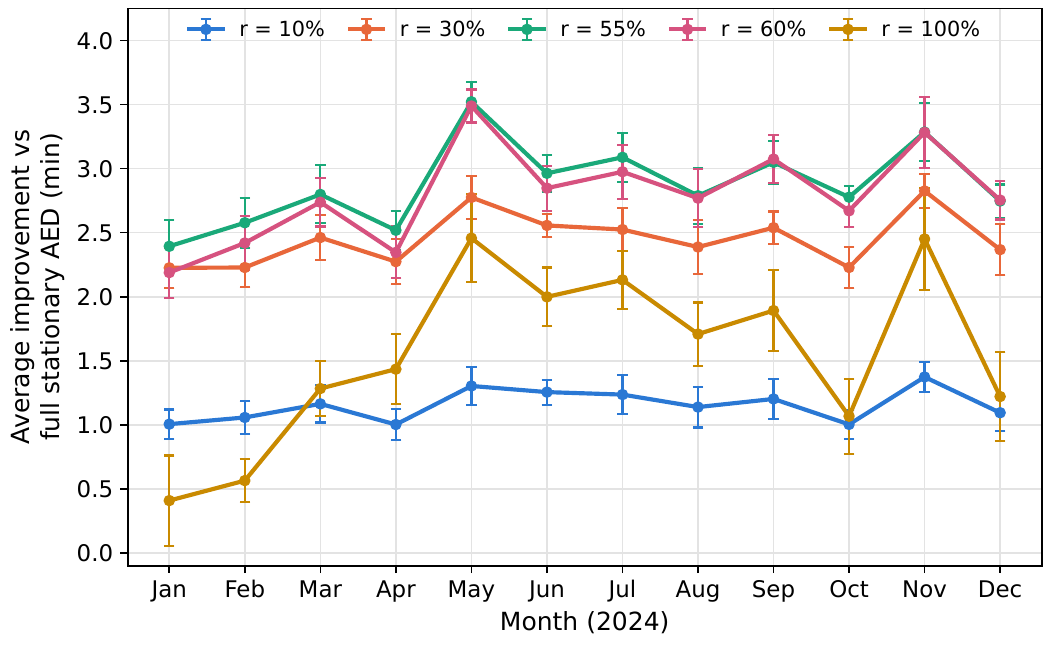}
    \par{\small (b) Mean improvement by hybrid ratios}
\end{minipage}
\caption{Monthly variation in mean response-time improvement (Toronto, 2024, occupied fleet). (a) Each line is one month, with the months of lowest (January) and highest (May) improvement highlighted; the dashed line is the 2024 average. (b) Monthly mean improvement at selected hybrid ratios; bars show $\pm$ one standard deviation across the ten mobility samples.}
\label{fig:toronto_monthly}
\end{figure}

In Toronto (Figure~\ref{fig:ratio_improvement}a), mean improvement in response time rises almost linearly at low hybrid ratios, reaching $2.45$~minutes at $30\%$, then increases more slowly to a maximum of $2.87$~minutes at $55\%$, and declines sharply beyond $60\%$ to $1.54$~minutes under a fully mobile policy. The optimum therefore keeps roughly half of the AED capacity at fixed locations.
Figure~\ref{fig:toronto_monthly} shows how these results vary across months and mobility samples. Panel (a) shows the monthly improvement curves which nearly coincide at low hybrid ratios and consistently peak between $50\%$ and $60\%$. Beyond this range, the monthly curves diverge but remain positive: further reallocation reduces the benefit without eliminating it. 
%Beyond this range, the monthly curves diverge, but improvements remain positive in every month, showing that further reallocation beyond the optimum reduces the benefit without eliminating the advantage over the full stationary policy.
Panel (b) provides a closer comparison at selected hybrid ratios. Intermediate allocations ($55\%$--$60\%$) consistently deliver great improvements throughout the year, with relatively small variation across the ten mobility samples within each month. In contrast, the fully mobile policy shows substantial monthly variation, with mean improvement ranging from $0.41$~minutes in January to $2.46$~minutes in May. Thus, an intermediate hybrid ratio performs well throughout the year, whereas the benefits of full mobilization are considerably more sensitive to monthly conditions. This partly reflects the seasonal variations of the mobility patterns.

% In New York City (Figure~\ref{fig:ratio_improvement}b), mean improvement rises steeply as AEDs are first moved to the mobile layer. At low hybrid ratios, the reallocated devices come from locations that retain another AED, so the mobile layer adds coverage without any loss of stationary access and its contribution is a pure gain, reaching $1.92$~minutes at $20\%$ and $2.61$~minutes at $40\%$.
% As reallocation begins to withdraw AED locations, the additional mobile coverage is increasingly offset by the stationary coverage it replaces, and mean improvement enters a broad plateau, remaining within about $10\%$ of its maximum between hybrid ratios of $40\%$ and $80\%$. The maximum, $2.91$ minutes, is reached at $80\%$.
% Beyond this point, the stationary losses dominate and mean improvement falls sharply, although a fully mobile network still improves on the stationary benchmark, by $1.14$~minutes.
% The hybrid policy therefore has an interior optimum that sits within a broad plateau. 
In New York City (Figure~\ref{fig:ratio_improvement}b), mean improvement reaches $1.92$~minutes at $20\%$ and $2.61$~minutes at $40\%$, then enters a broad plateau, remaining within about $10\%$ of its maximum between hybrid ratios of $40\%$ and $80\%$. The maximum, $2.91$~minutes, is reached at $80\%$; beyond it, mean improvement falls sharply, to $1.14$~minutes under a fully mobile policy. The optimum therefore sits within a broad plateau.
This pattern is stable across the year, with mean improvement peaking at a hybrid ratio of $75\%$ or $80\%$ in every month of 2024, and the monthly peak ranging from $2.76$ to $3.29$~minutes.
%The breadth of the plateau gives considerable flexibility. Because mean improvement changes little across a wide range of hybrid ratios, the choice within that range can be guided by other objectives, such as limiting the number of OHCA events made worse off, as discussed next.

The two cities therefore attain a similar maximum mean improvement in response time, 2.87 minutes in Toronto and 2.91 minutes in New York City, despite very different stationary coverage.
The key difference lies in the flexibility of deployment decisions: Toronto's relatively narrow optimum calls for a more carefully balanced stationary--mobile allocation, whereas New York City's broad plateau provides the planner with great flexibility. For example, a planner may favor the lower end of the plateau to preserve stationary coverage while retaining most of the attainable response-time improvement.

\subsubsection{Distribution of Gains and Losses}

\begin{figure}[!b]
\centering
\begin{minipage}[t]{0.49\textwidth}
    \centering
    \includegraphics[width=\textwidth]{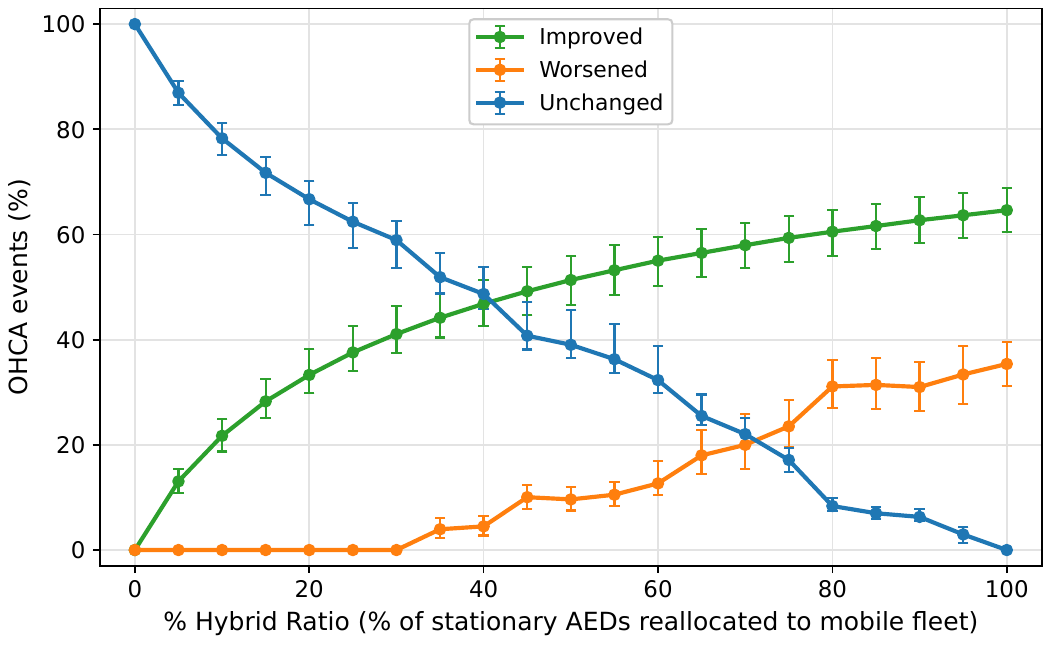}
    \par\smallskip{\small (a) Toronto, 2024 }
\end{minipage}\hfill
\begin{minipage}[t]{0.49\textwidth}
    \centering
    \includegraphics[width=\textwidth]{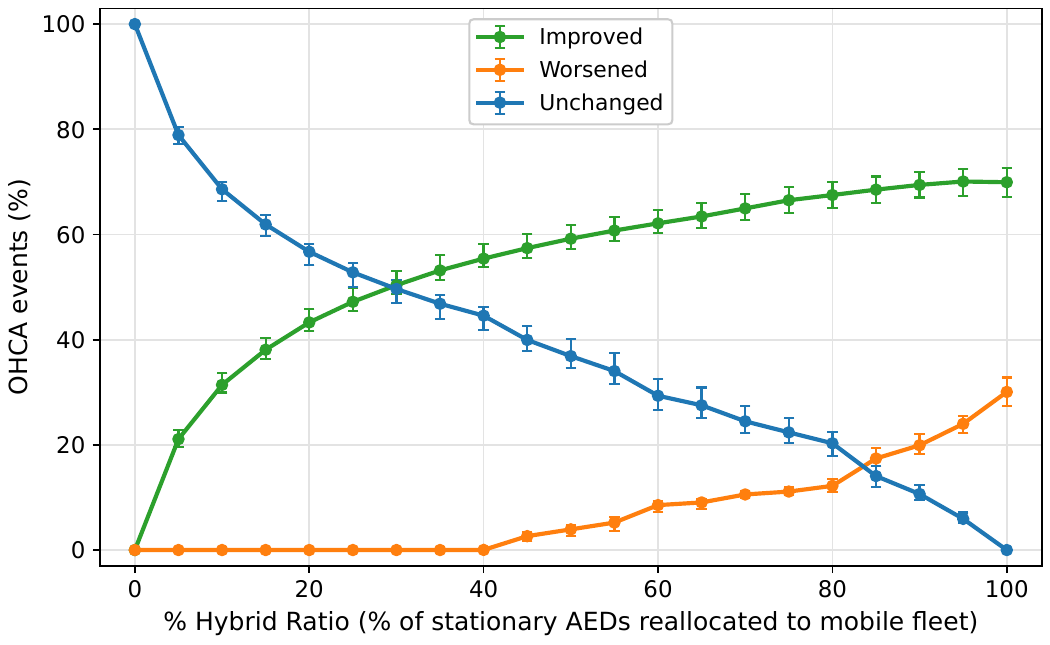}
    \par\smallskip{\small (b) New York City, 2024}
\end{minipage}
\caption{Share of OHCA events with improved, worsened, and unchanged response time as a function of the hybrid ratio, under the occupied dispatch rule (2024): (a) Toronto, (b) New York City. Bars show the range across months.}
\label{fig:ratio_counts}
\end{figure}
%. Events with no mobile or stationary response are not shown; they account for less than 0.01\% of events at every hybrid ratio in both cities

Another policy-relevant question is how broadly the benefit is distributed. 
We therefore adopt a binary improvement measure and classify each OHCA event as a gain (improved response) or a loss (worsened response) relative to the full stationary network, and report the share of events in each class across the range of hybrid ratios, averaged across mobility samples. (Figure~\ref{fig:ratio_counts}).

The share of ``gain'' events rises at every hybrid ratio in both cities, steeply as the first AEDs are moved to vehicles and more slowly thereafter, until roughly two-thirds of events are reached faster than under the full stationary network once all AEDs are mobile. 
Additional equipped vehicles therefore keep extending the benefit to new events, but at a diminishing rate.

Losses follow a different course. They are absent at low hybrid ratios and appear only once reallocation begins to withdraw non-duplicate AED locations; the share of worsened events then rises in steps, and more steeply at high hybrid ratios, as events that a nearby fixed device would otherwise have served lose that device. 

% The loss-free range, up to $30\%$ in Toronto and $35\%$ in New York City, follows from the reallocation order, since the additional devices at locations hosting more than one are reallocated first, no location is withdrawn until they are exhausted.

In Toronto (Figure~\ref{fig:ratio_counts}a), the improved share grows from $41.1\%$ at a hybrid ratio of $30\%$ to $64.6\%$ under a fully mobile policy. The worsened share reaches $10.5\%$ at $55\%$ and $35.4\%$ under a fully mobile policy.
In New York City (Figure~\ref{fig:ratio_counts}b), the improved share grows from $43.3\%$ at $20\%$ to $69.9\%$ under a fully mobile policy. 
The worsened share reaches $5.2\%$ at $55\%$ and $30.1\%$ under a fully mobile policy.

Mean improvement therefore falls beyond the optimum even though the share of gains keeps rising, because losses grow in magnitude faster than gains. In Toronto, events that gain at a hybrid ratio of $55\%$ receive an AED $6.5$~minutes sooner on average, and $7.1$~minutes sooner under a fully mobile policy, whereas events that lose receive it $5.5$~minutes later at $55\%$ and $8.6$~minutes later under a fully mobile policy. Between these two ratios, the average loss grows by $3.1$~minutes while the average gain grows by only $0.6$~minutes.

Read together with the plateau in mean improvement, this composition shows a trade-off in both cities between the size of the average benefit and the number of events made worse off. Moving from a hybrid ratio of $40\%$ to the ratio that maximizes mean improvement raises mean improvement by $0.26$~minutes in Toronto (to $55\%$) and $0.30$~minutes in New York City (to $80\%$), while the worsened share rises from $4.5\%$ to $10.5\%$ in Toronto and from almost none to about one in eight in New York City. A planner who weighs the protection of events already well served by the stationary network can therefore choose a hybrid ratio near $40\%$ in either city and retain about $90\%$ of the attainable benefit.

\subsubsection{Reliability}

\begin{figure}[!b]
\centering
\begin{minipage}[t]{0.49\textwidth}
    \centering
    \includegraphics[width=\textwidth]{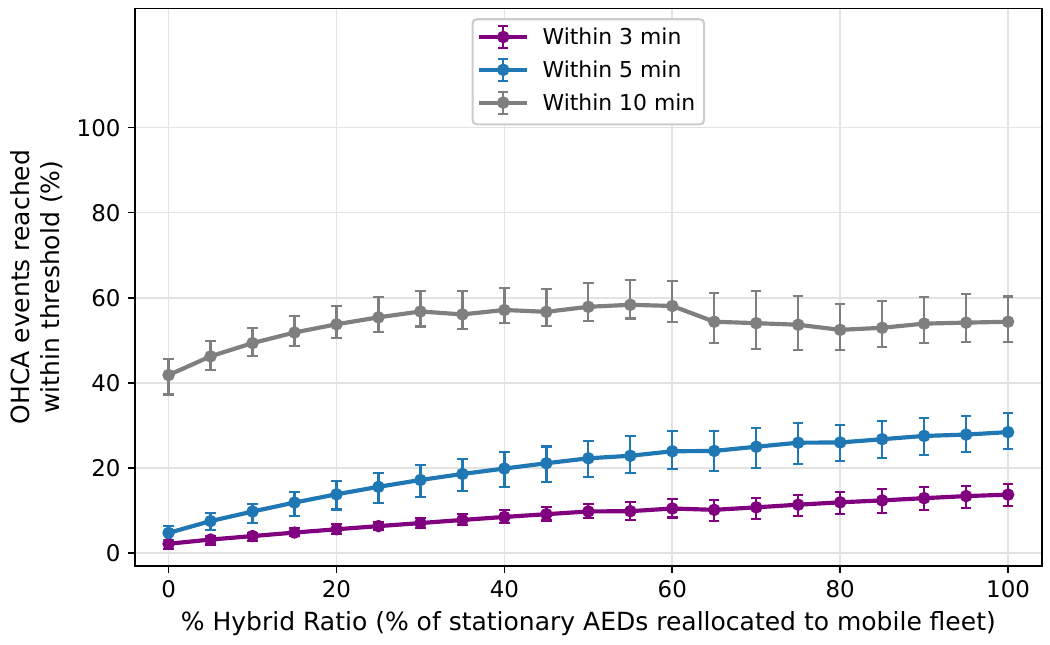}
    \par\smallskip{\small (a) Toronto, 2024}
\end{minipage}\hfill
\begin{minipage}[t]{0.49\textwidth}
    \centering
    \includegraphics[width=\textwidth]{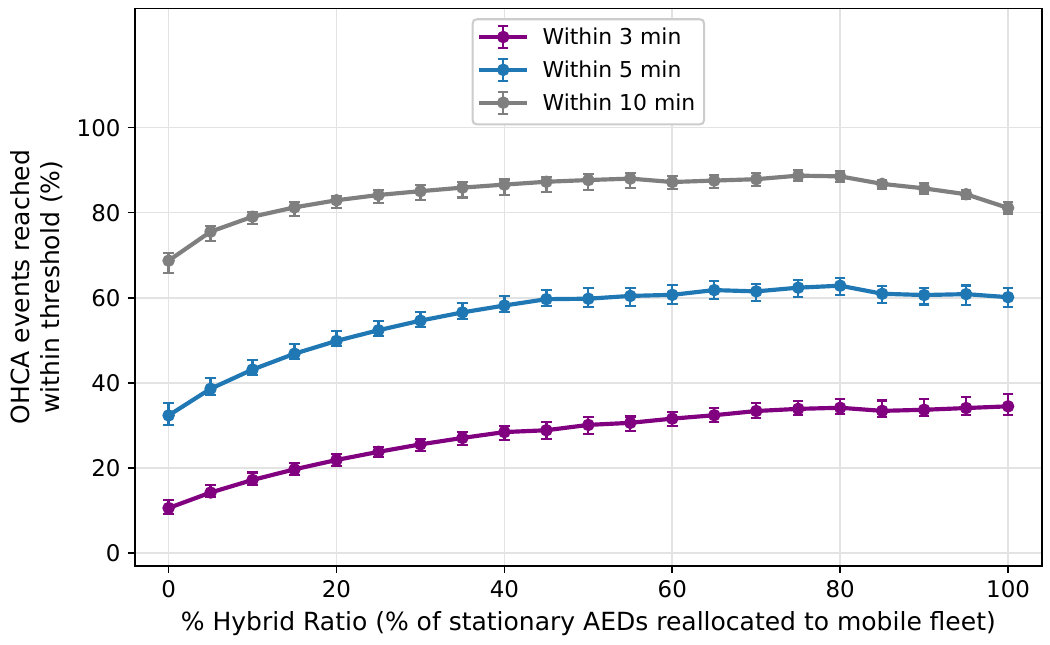}
    \par\smallskip{\small (b) New York City, 2024}
\end{minipage}
\caption{Reliability, share of OHCA events reached within $3$, $5$, and $10$~minutes, as a function of the hybrid ratio, under the occupied dispatch rule (2024): (a) Toronto, (b) New York City. Bars show the range across months.}
\label{fig:reliability}
\end{figure}

Reliability, defined as the share of OHCA events reached within 3, 5, and 10 minutes, exhibits different patterns across response-time thresholds (Figure~\ref{fig:reliability}). In both cities, reliability within 10 minutes peaks at an intermediate hybrid ratio, whereas shorter response-time thresholds continue to benefit from greater mobile deployment. This difference arises because short response times are achieved primarily by mobile AEDs, while 10-minute coverage also depends on stationary AEDs serving events close to a fixed device, which the mobile layer does not fully replace once those locations are withdrawn.
% \textcolor{gray}{In Toronto (Figure~\ref{fig:reliability}a), the share of events reached within 10 minutes rises from 41.8\% under the stationary policy to a maximum of 58.3\% at a hybrid ratio of 55\%, then declines at higher ratios. In contrast, the share reached within 5 minutes continues to increase, reaching almost five times its stationary level at 55\%.
% New York City exhibits a similar pattern (Figure~\ref{fig:reliability}b). The share of events reached within 10 minutes rises from 68.7\% to almost 89\%, then declines at high hybrid ratios, whereas reliability at shorter thresholds remains at or near its maximum up to a fully mobile network. Unlike Toronto, reliability at every threshold remains above the full stationary baseline across all hybrid ratios.}
The preferred hybrid ratio therefore depends on the response-time target: achieving the shortest response times favors greater mobile deployment, whereas maintaining broad coverage within 10 minutes favors retaining part of the stationary network.

\subsubsection{Heterogeneity Within New York City}

The borough-level analysis in Appendix~\ref{app:nyc-borough} shows that the citywide curve is a composite of distinct local responses. The hybrid ratio that maximizes mean improvement varies from 40\% in Richmond/Staten Island, where ride-hailing activity is lowest, to 85\% in the Bronx, and under a fully mobile policy, access improves in four boroughs but deteriorates substantially in Richmond/Staten Island. The largest gains occur in boroughs that combine moderate stationary coverage with substantial fleet activity, such as Brooklyn ($3.9$~minutes at its optimum), rather than where coverage is weakest (Richmond/Staten Island) or strongest (Manhattan). The appropriate balance between stationary and mobile supply therefore varies within the city, and a single citywide hybrid ratio can place parts of the city in a regime where further reallocation reduces access.

\subsection{Policy Implications: Dispatch and Allocation in Toronto}
\label{sec:policy-implications}

The Toronto records identify all three activity phases and a vehicle identifier, allowing us to examine two operational levers: which vehicles may be dispatched and which vehicles are equipped.

\subsubsection{Dispatch Policy Across Ride-Hailing Trip Phases}
\label{sec:toronto-availability} 

%\paragraph{Idle versus occupied vehicles.}]
\begin{figure}[!b]
    \centering
    \includegraphics[width=0.8\linewidth]{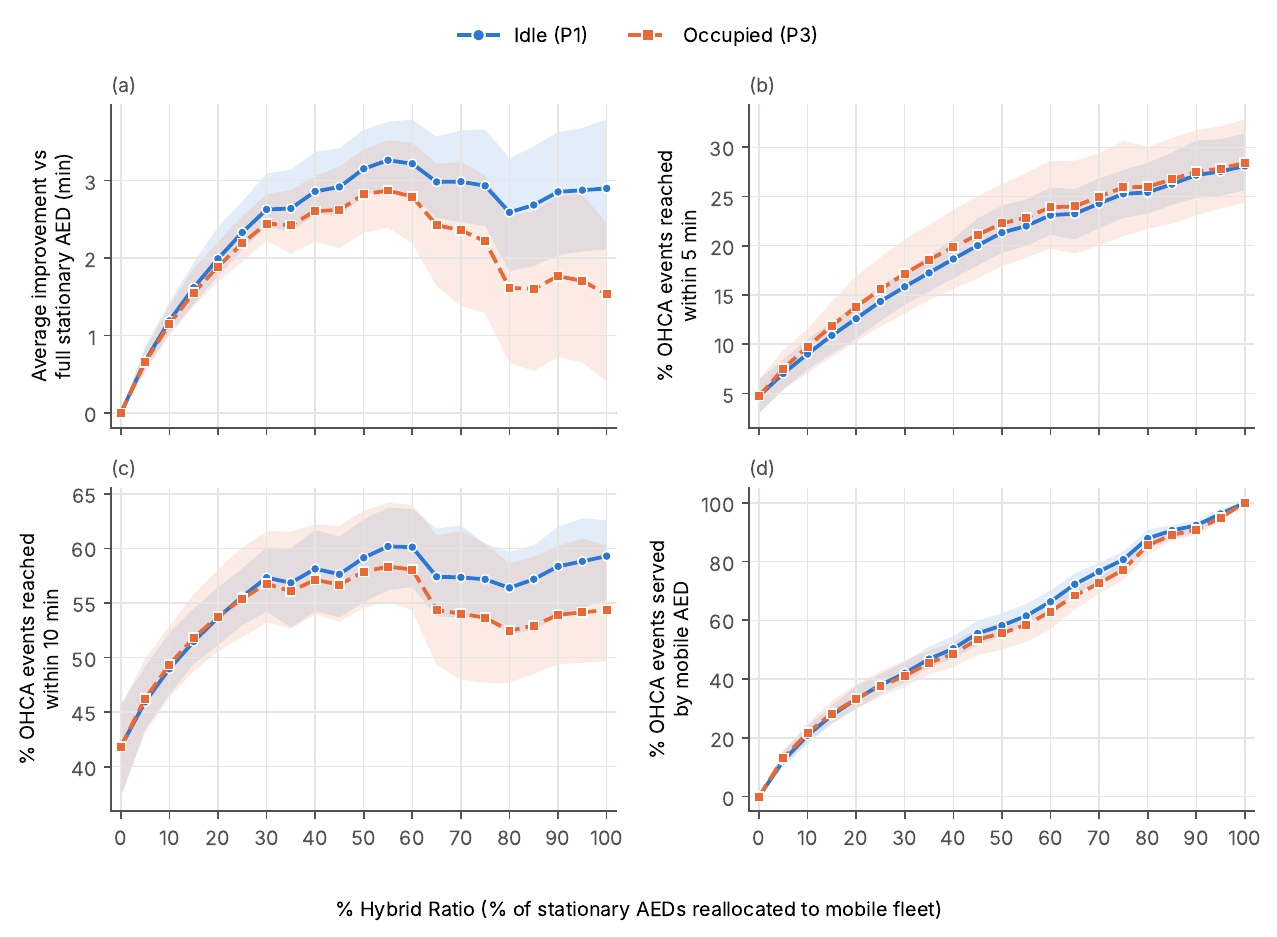}
    \caption{Idle (P1) versus occupied (P3) dispatch rule as a function of the hybrid ratio (Toronto, 2024): (a) mean response-time improvement relative to the full stationary network; share of OHCA events reached within (b) $5$ and (c) $10$~minutes; (d) share of OHCA events served by a mobile AED. Lines are averages over the year and the ten mobility samples; shaded bands show the range across months. The equipped vehicles are identical under both rules.}
    \label{fig:toronto_rules_p3_p1}
\end{figure}

An equipped vehicle carries its AED through all three activity phases, so the phases do not change which vehicles hold AEDs, only which of them may respond at the instant of an OHCA. We refer to the set of phases in which an equipped vehicle may be dispatched as the \emph{dispatch rule}. In addition to the occupied rule of Section~\ref{sec:toronto-hybrid}, we consider the \emph{idle} rule, under which only equipped vehicles in P1 respond, so that no passenger trip is disrupted, and the \emph{all-phase} rule, under which any equipped vehicle responds.

Figure~\ref{fig:toronto_rules_p3_p1} compares the idle rule with the occupied rule of the main analysis. At low hybrid ratios the two rules produce similar outcomes. Although the two rules almost never dispatch the same vehicle, mean improvement differs by at most $0.1$~minutes up to a hybrid ratio of $20\%$. The idle rule then pulls ahead. At $55\%$, which is the peak for both rules, the mean response improvement is $3.26$~minutes under the idle rule against $2.87$~minutes under the occupied rule, and the difference widens at higher hybrid ratios.
%, reaching $2.90$ against $1.54$~minutes under a fully mobile policy. 
The idle fleet also reaches more OHCA events at the $10$~minutes reliability threshold, and its advantage grows at high hybrid ratios, where the mobile fleet must compensate for withdrawn stationary locations. %, $60.2\%$ against $58.3\%$ at $55\%$ and $59.3\%$ against $54.4\%$ under a fully mobile policy. 
At the $5$-minute reliability threshold, the ordering reverses slightly, and the occupied fleet is marginally ahead. 
%Under activity-weighted equipping,
Therefore, limiting dispatch to vehicles that can be diverted without interrupting a passenger trip does not lower mean improvement relative to the occupied specification of the main analysis, and on average that specification does not overstate what a non-disruptive policy can achieve. This comparison is nevertheless sensitive to how vehicles are equipped, as discussed in Section~\ref{sec:toronto-selection} and ~\ref{app:toronto-data}.

\begin{figure}[!b]
    \centering
    \includegraphics[width=0.8\linewidth]{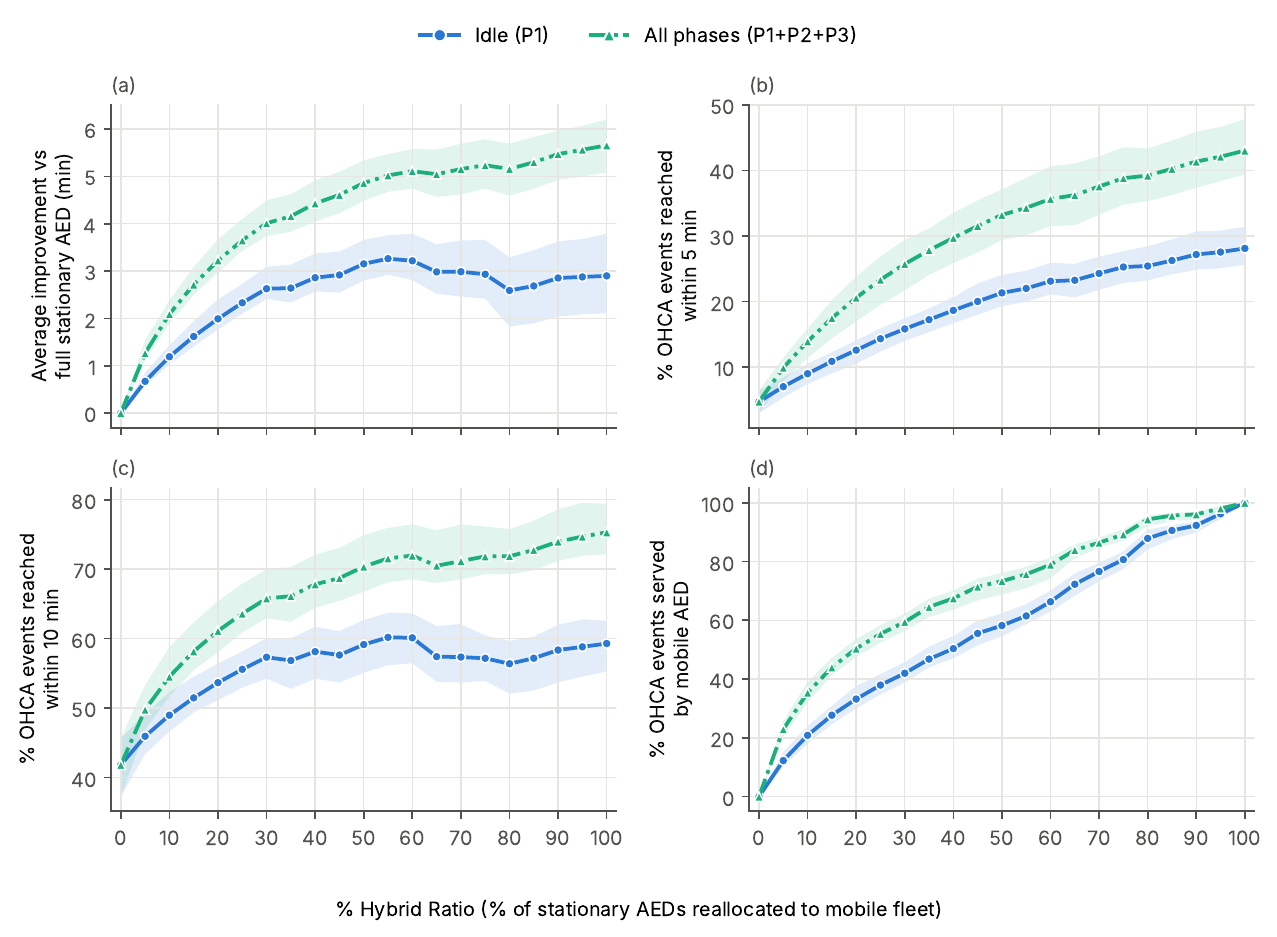}
    \caption{Idle (P1) versus all-phase (P1+P2+P3) dispatch rule as a function of the hybrid ratio (Toronto, 2024): (a) mean response-time improvement relative to the full stationary network; share of OHCA events reached within (b) $5$ and (c) $10$~minutes; (d) share of OHCA events served by a mobile AED. Lines are averages over the year and the ten mobility samples; shaded bands show the range across months. The equipped vehicles are identical under both rules.}
    \label{fig:toronto_rules_all_p1}
\end{figure}

%\paragraph{All phases versus idle vehicles.}

The idle rule is a restriction of the all-phase rule, and the gap between them measures the cost of not dispatching equipped vehicles en route to a pickup or carrying a passenger (Figure~\ref{fig:toronto_rules_all_p1}). This gap is large at every hybrid ratio and on every measure, and it stems from fleet size: roughly three times as many equipped vehicles are available at the instant of an OHCA under the all-phase rule. 
As a result, the all-phase rule at a hybrid ratio of $25\%$ already exceeds the largest mean improvement the idle rule can achieve. The all-phase rule also changes the shape of the response to reallocation, and the mean improvement no longer peaks at an intermediate ratio but continues to increase, apart from small reversals, up to a fully mobile policy.
% \textcolor{gray}{Under the idle and occupied rules mean improvement has an interior maximum at $55\%$, whereas under the all-phase rule it continues to increase, with small reversals at $65\%$ and $80\%$, to $5.66$~minutes under a fully mobile policy, and the $10$-minute share rises to $75.3\%$. Mobile AEDs serve $75.7\%$ of OHCA events at $55\%$, against $61.5\%$ under the idle rule, and the worsened share is roughly halved, from $11.2\%$ to $5.9\%$ at $55\%$ and from $31.7\%$ to $18.8\%$ under a fully mobile policy.}
However, a fully mobile policy with all-phase dispatching may still leave a share of OHCA events reached later than under the full stationary network, so the choice of hybrid ratio remains a trade-off between average improvement and worse-off events. From a policy perspective, the gap between the two curves is the value of permitting the diversion of vehicles carrying, or about to collect, a passenger, and bounds what a platform could offer drivers and passengers in compensation for such diversions.
Thus, the interior optimum appears whenever the available equipped fleet is too small to replace the stationary AEDs withdrawn at high hybrid ratios. %A restricted dispatch rule and uniform equipping both produce this, while the all-phase rule with activity-weighted equipping does not.

\subsubsection{Which Vehicles To Be Equipped}
%Sensitivity to Fleet Availability
\label{sec:toronto-selection}

In the analysis above, AEDs are assigned to vehicles in proportion to their activity on the platform, which %a city can implement in partnership with ride-hailing platforms that hold this information and 
can shape through driver incentives. To assess the value of this targeting, we repeat the analysis with the equipped vehicles drawn uniformly at random from the vehicles with at least one occupied trip in the same month, keeping the counterfactual AED network and all other elements unchanged. The difference between the two specifications therefore isolates the effect of targeted vehicle distribution.

\begin{figure}[!b]
    \centering
    \includegraphics[width=0.8\linewidth]{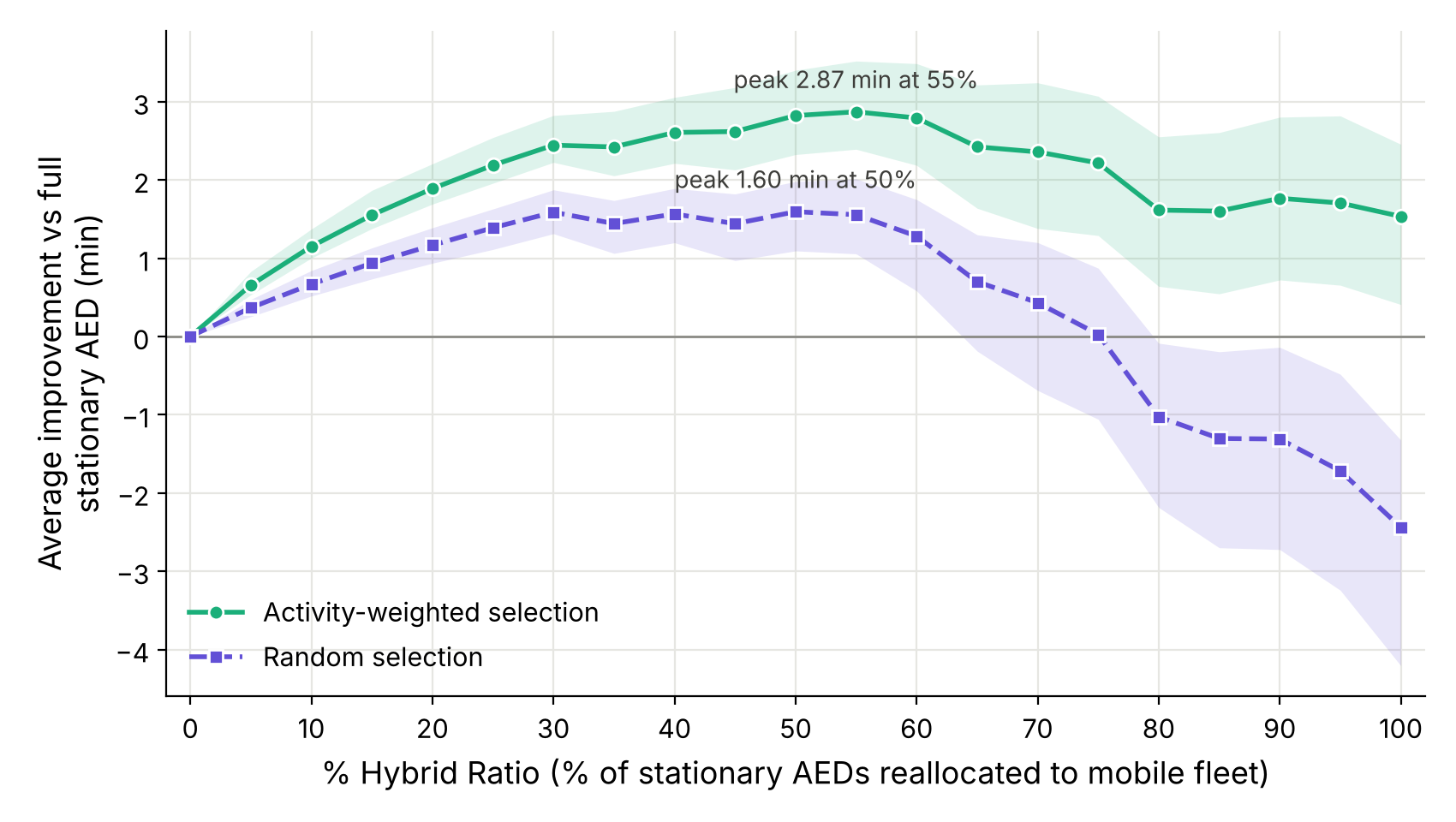}
    \caption{Mean response-time improvement relative to the full stationary network under activity-weighted and uniform equipping of equipped vehicles (Toronto, 2024, occupied fleet). Bands show the range across months.}
    \label{fig:toronto_selection}
\end{figure}

% Vehicle selection changes how the benefit evolves with the hybrid ratio (Figure~\ref{fig:toronto_selection}). Under uniform equipping (random selection), mean improvement reaches approximately $1.6$~minutes at a hybrid ratio of $30\%$ and does not increase further, compared with a peak of $2.87$~minutes at $55\%$ under activity-weighted selection. Beyond a hybrid ratio of about $75\%$, uniform equipping performs worse than the full stationary network, with mean response time approximately $2.4$~minutes slower under a fully mobile policy, whereas activity-weighted selection still improves on the stationary network by $1.5$~minutes. Reliability follows the same pattern: under uniform equipping, the share of events reached within $10$~minutes returns to its stationary level of about $42\%$ under a fully mobile policy. 

Figure~\ref{fig:toronto_selection} shows how vehicle selection affects the benefits of mobile AED deployment. Under uniform equipping, mean improvement quickly levels off and eventually becomes negative as more stationary AEDs are reallocated, whereas activity-weighted selection achieves greater improvements and maintains positive gains even under a fully mobile policy. Reliability follows a similar pattern, highlighting the importance of targeting more active vehicles when expanding mobile AED deployment.

A mobile AED contributes only when its vehicle is on the road at the time of an arrest, while uniform equipping places many devices in vehicles that are rarely active.
Targeted distribution helps ensure the AEDs are allocated to vehicles more likely to be present close to an OHCA. Without vehicle-activity information, mobile AED gains saturate at a low hybrid ratio, and further reallocation becomes detrimental, so a larger share of AEDs should remain at fixed locations. Therefore, the city should carefully design the allocation rule and leverage activity information from the ride-hailing platform to maximize the gains from mobile AEDs. %In the presence of vehicle activity information, the city can move more capacity to the mobile layer and obtain a larger benefit at lower risk.

 \subsection{Vehicle Equipping Sensitivity Analysis}
 \label{app:toronto-data}
\begin{figure}[b]
    \centering
    \includegraphics[width=\linewidth]{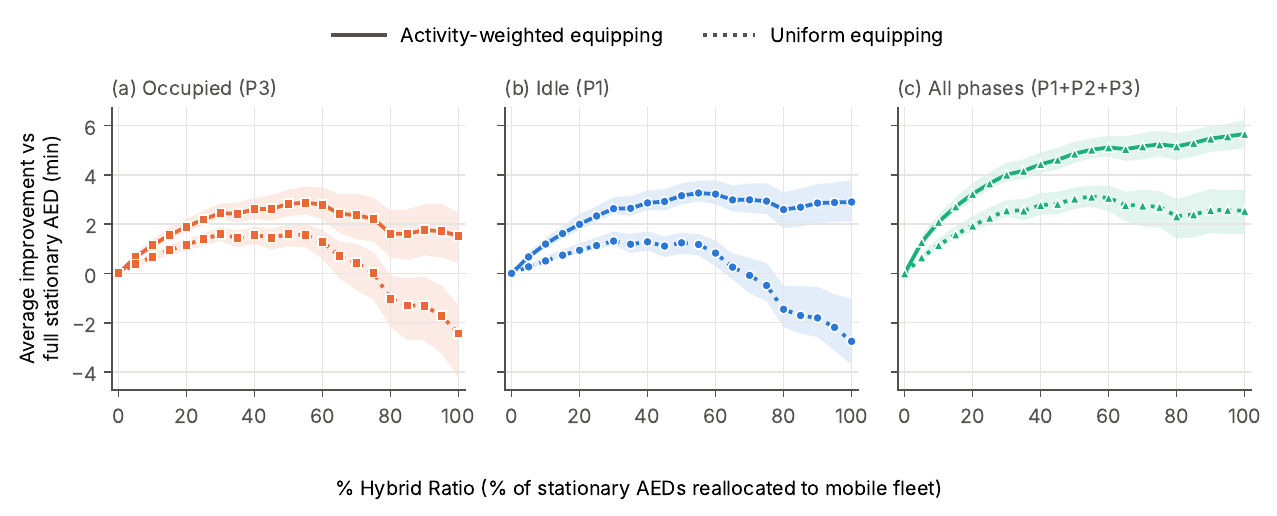}
    \caption{Mean response-time improvement under activity-weighted (solid) and uniform (dotted) equipping, by dispatch rule (Toronto, 2024): (a) occupied (P3), (b) idle (P1), (c) all-phase (P1+P2+P3). Lines are averages over the year and the ten mobility samples; shaded bands show the range across months.}
    \label{fig:toronto_equipping}
\end{figure}

The results of Section~\ref{sec:toronto-availability}, like those of Section~\ref{sec:toronto-hybrid}, draw the equipped vehicles with probability proportional to their activity in the month, from the vehicles with at least five occupied trips in that month. Figure~\ref{fig:toronto_equipping} repeats the analysis with random equipping, in which every vehicle with at least one occupied trip in the month is equally likely to be equipped. Uniform equipping can place a larger share of the AEDs on vehicles that are rarely active, so fewer equipped vehicles are available at the instant of an OHCA, and performance falls under every rule. %Peak mean improvement falls to $1.60$~minutes at $50\%$ under the occupied rule, $1.31$~minutes at $30\%$ under the idle rule, and $3.12$~minutes at $55\%$ under the all-phase rule. 
Under the two restricted rules, mean improvement becomes negative at high hybrid ratios, whereas under the all-phase rule it remains positive at every ratio, does not fall below $2.3$~minutes beyond its peak, and now exhibits an interior optimum.

Two comparisons change under uniform equipping. First, the ordering of the idle and occupied rules reverses: with uniform equipping the occupied rule yields the larger mean improvement at every hybrid ratio. The relative performance of the two restricted rules therefore depends on which vehicles carry the AEDs, whereas the advantage of the all-phase rule holds under both equipping rules. Second, the all-phase rule with uniform equipping performs close to the idle rule with activity-weighted equipping ($3.12$ against $3.26$~minutes at $55\%$). 
The choice of which vehicles to equip is therefore a policy lever comparable to the choice of dispatch rule. Taken together, the dispatch and equipping comparisons indicate when the interior optimum arises: it appears whenever the available equipped fleet is too small to replace the stationary AEDs withdrawn at high hybrid ratios. A restricted dispatch rule and uniform equipping both produce this, while the all-phase rule with activity-weighted equipping does not.
%an idle-only rule with activity-weighted equipping attains what the all-phase rule attains with uniform equipping, while the all-phase rule with activity-weighted equipping remains well above both.

% Both changes act on the number of equipped vehicles available when an OHCA occurs, one by admitting more phases and the other by placing AEDs on more active vehicles, which suggests that this number, rather than the phase or the equipping rule as such, is the main driver of the benefit. % verify mean available equipped vehicles per event for all six runs; if mean improvement lines up on one curve against it, this sentence stands

\section{Methods}
\label{sec:methods}

The analysis was implemented in Python, using pandas, NumPy, GeoPandas, and DuckDB for data processing. All road-network distances were computed with the Open Source Routing Machine (OSRM)~\citep{luxen2011osrm} on OpenStreetMap road data~\citep{osm}. We analyze both cities over the calendar year 2024 using a common construction; when the data for the two cities support different treatments, we state each separately. Throughout, $\mathcal{E}$ denotes the set of OHCA events, and event $e\in\mathcal{E}$ occurs at location $x_e$ and time $t_e$.

\subsection{Events, AEDs, and Vehicle Positions}
\label{sec:methods-data}

The construction of the three data layers is described in Appendix~\ref{app:data}. Each OHCA event is placed at a representative point of its postal area, the population-weighted centroid of its FSA in Toronto and the centroid of its ZIP code in New York City, and each AED at its geocoded location. We write $D$ for the total number of AED devices in a city's inventory, the capacity held fixed across all policies ($D = 1{,}291$ in Toronto and $11{,}138$ in New York City). Each Toronto segment and each New York City trip is routed on the driving network between its endpoints, and a vehicle's position at any instant is interpolated along its routed path at the observed pace of the segment or trip. For a segment $i$ that starts at time $a_i$, ends at time $b_i$, and follows a routed path $\gamma_i$ of length $L_i$, parameterized by distance traveled, the position at time $t\in[a_i,b_i]$ is
\begin{equation}
x_i(t) \;=\; \gamma_i\!\left(L_i\,\frac{t-a_i}{b_i-a_i}\right).
\label{eq:interp}
\end{equation}
Total travel time therefore matches the observed duration, and only the pace within a segment is approximated as constant.

\subsection{Routing and Travel Time Estimation}
\label{sec:methods-routing}

We run separate OSRM instances for the driving and foot profiles in each city. For Toronto, the road network is extracted for a bounding box of the Greater Toronto Area rather than the city boundary, so that trips that start or end outside the city, for example at Toronto Pearson International Airport, are routed rather than snapped to the boundary. For New York City, the road network is the New York State OpenStreetMap extract.
Each instance was verified on known origin--destination pairs (implied speed, one-way asymmetry, turn restrictions, and pedestrian detour) before use. We write $d_F(x,y)$ for the on-foot network distance and $d_D(x,y)$ for the driving network distance from $x$ to $y$.

\paragraph{Stationary leg.}
Let $\mathcal{K}_r$ denote the set of AED locations that retain at least one device at hybrid ratio $r$ (defined in the Allocation Rule below). The stationary leg of event $e$ is a two-way retrieval on foot to the nearest retained AED:
\begin{equation}
d^{S}_{r}(e) = \min_{a\in\mathcal{K}_r} d_F(x_e, a),
\qquad
\tau^{S}_{r}(e) = \frac{2\,d^{S}_{r}(e)}{60\,v_F},
%time unit?!
\label{eq:stat}
\end{equation}
with walking speed $v_F = 1.3$~m/s. % Bahar: add the walking-speed reference (same Manhattan source).

\paragraph{Mobile leg.}
Let $\mathcal{A}_{r,v}(t)$ denote the set of equipped vehicles that are available to respond at time $t$ at hybrid ratio $r$ in mobility sample $v$ (defined below). The mobile leg is a one-way drive from the nearest available equipped vehicle to the event:
\begin{equation}
d^{M}_{r,v}(e) = \min_{i\in\mathcal{A}_{r,v}(t_e)} d_D\big(x_i(t_e),\, x_e\big),
\qquad
\tau^{M}_{r,v}(e) = \frac{d^{M}_{r,v}(e)}{60\,v_D},
\label{eq:mob}
\end{equation}
with vehicle speed $v_D = 3.17$~m/s ($11.4$~km/h), the 2017 average vehicle travel speed in the Manhattan central business district~\citep{NYCDOT2018}, which we use as a conservative vehicle speed in both cities. % Bahar: add the Manhattan speed reference.
If no equipped vehicle is available at $t_e$, we set $\tau^{M}_{r,v}(e) = \infty$. We use this constant speed rather than OSRM's travel-time estimates, which assume free-flow conditions and imply a median vehicle speed of about $14$~m/s in Toronto. Combined with the assumptions that a second bystander is always available and that a stationary AED is retrieved without delay for locating or unlocking it, this choice favours the stationary benchmark.
% two other assumptions: The mobile leg carries no notification, acceptance, or parking delay, and every registered AED is assumed reachable at all hours (the latter favours the stationary side). A per-mode fixed activation offset is the natural sensitivity if a reviewer raises the first.

%\paragraph{Computation.}
In Toronto, the stationary leg is computed from the complete on-foot distance matrix between the $96$ event locations and all $902$ AED locations, so the minimum in Equation~\eqref{eq:stat} is exact at every ratio. For the mobile leg, every vehicle active on the platform at $t_e$, in any phase, is routed to the event location through the OSRM table service, without pre-selection. These distances are computed once and reused for all hybrid ratios, mobility samples, and dispatch rules, which differ only in which routed vehicles are admitted to the minimum in Equation~\eqref{eq:mob}. %New York City follows the same construction at a larger scale, with the nearest retained AED and the nearest available equipped vehicle pre-selected by straight-line distance before being measured on the network.
% \textcolor{gray}{
% Because the retained and equipped sets are nested across ratios (see below), a distance measured at one ratio remains attainable at every ratio whose set contains the measured candidate. We therefore replace each measured distance by the smallest measured distance among those ratios,
% \begin{equation}
% \tilde{d}^{S}_{r}(e) = \min_{r'\ge r} d^{S}_{r'}(e),
% \qquad
% \tilde{d}^{M}_{r,v}(e) = \min_{r'\le r} d^{M}_{r',v}(e),
% \label{eq:monotone}
% \end{equation}
% which uses only routed distances and ensures that the stationary leg cannot shorten and the mobile leg cannot lengthen as the hybrid ratio increases.}
% Bahar: confirm that the 12-month NYC outputs reported in the paper include this correction (overview notes it as adopted). If they do not, delete this paragraph, but the NYC results would then need a check for rows where the rolled-back stationary leg beats the full-inventory baseline.

\subsection{Allocation Rule}
\label{sec:methods-allocation}

\paragraph{Hybrid ratio.}
The hybrid ratio $r\in\{0, 5, 10, \dots, 100\}$ is the percentage of the $D$ AEDs allocated to the mobile layer. At ratio $r$,
% \begin{equation}
% n_r = \operatorname{round}\!\left(\frac{r}{100}\,D\right)
% \label{eq:nr}
% \end{equation}
the number of devices that are reallocated from the stationary network to equipped vehicles is $n_r = \operatorname{round}\!\left(\frac{r}{100}\,D\right)$, and $D - n_r$ devices remain at fixed locations, so total AED capacity is identical under every policy. Both the stationary rollback and the vehicle equipping are each defined by a single ordering that is cut at $n_r$. As a result, the configurations are nested: for $r < r'$, every location retained at $r'$ is retained at $r$, and every vehicle equipped at $r$ is equipped at $r'$. %A reallocated device is not removed from the system; in the counterfactual, it is carried by a vehicle.

\paragraph{Stationary rollback.}
We construct the counterfactual stationary networks by reversing densification, the addition of devices at locations that already host one, and apply the same rule in both cities. Devices are withdrawn in two stages. First, the additional devices at locations hosting more than one device are withdrawn, leaving one device at each location. Second, each location's remaining device is withdrawn, with locations of protected types withdrawn last: schools and colleges, health and social services including long-term care, government facilities including police, fire, and paramedic services, and transit stations. Within each stage, the order is set by a random priority score, drawn with a fixed seed, that favors devices in postal areas (FSAs in Toronto, ZIP codes in New York City) with more AED locations, so that reallocation starts where stationary coverage is densest. Under this order, no AED location is withdrawn while $n_r$ does not exceed the number of co-located devices, that is, up to a hybrid ratio of $30.1\%$ in Toronto and $33.8\%$ in New York City. The rollback is deterministic given the seed, and one stationary configuration per ratio is shared by all mobility samples. %bahar add reference to tim's paper for relocating first from redundant locations.

\paragraph{Vehicle equipping.}
In Toronto, mobile AEDs are assigned to vehicles (with anonymized vehicle IDs). For month $m$, the eligible pool consists of the vehicles with at least five completed occupied trips in that month. In each mobility sample $v$, the pool is put in a random order by weighted sampling without replacement: at each draw, a remaining vehicle $i$ is selected with probability

\begin{equation}
\pi_i = \frac{w_{i}}{\sum_{j\ \text{remaining}} w_{j}},
\label{eq:weighted}
\end{equation}
where $w_i$ is the time vehicle $i$ spends idle on the platform in month $m$, so that AEDs are placed preferentially on the vehicles most often available to respond. The equipped set at ratio $r$ consists of the first $n_r$ vehicles in this order. Under the uniform equipping used for comparison, $w_i = 1$ for every vehicle with at least one occupied trip in the month. 
%In the hybrid_rule.py, w_i was the vehicle's idle (P1) hours in the month, applied to the same fleet for all three dispatch rules. we can write "the time vehicle i spends idle on the platform in month m". "Active" is correct when switching to total online hours.

In New York City, the trip records carry no vehicle identifier, so mobile AEDs are assigned to trips. The $n_r$ equipped vehicles represent a share
\begin{equation}
f_{r,m} = \frac{n_r}{V_m}
\label{eq:share}
\end{equation}
of the ride-hailing fleet, where $V_m$ is the average number of vehicles active per day in month $m$ reported by the TLC; $f_{100,m}$ ranges from $18.3\%$ to $20.2\%$ across the months of 2024. In each mobility sample, the month's trips are put in a random order and the first share $f_{r,m}$ of them are equipped. Because each trip is equally likely to be equipped, the mobile layer follows occupied ride-hailing activity in proportion to its volume. In both cities, ten mobility samples $v = 1,\dots,10$ are drawn per month with independent fixed seeds.

\paragraph{Dispatch rule.}
A dispatch rule specifies the set $\Phi$ of activity phases in which an equipped vehicle may respond, so that
\begin{equation}
\mathcal{A}_{r,v}(t) = \big\{\, i \in \mathcal{V}_{r,v} \,:\, \text{vehicle } i \text{ is active at } t \text{ in a phase in } \Phi \,\big\},
\label{eq:avail}
\end{equation}
where $\mathcal{V}_{r,v}$ is the equipped set. The occupied rule uses $\Phi = \{\mathrm{P3}\}$, the idle rule $\Phi = \{\mathrm{P1}\}$, and the all-phase rule $\Phi = \{\mathrm{P1}, \mathrm{P2}, \mathrm{P3}\}$. In Toronto, the equipped set does not depend on the rule, so the three rules differ only in which equipped vehicles may respond. In New York City, only the occupied phase is observed, so an equipped trip is available at $t$ if $a_i \le t \le b_i$.

\paragraph{Response rule.}
Under every policy, the AED that reaches the event first responds. The response time of event $e$ is
\begin{equation}
T_{r,v}(e) = \min\big\{\tau^{S}_{r}(e),\; \tau^{M}_{r,v}(e)\big\},
\label{eq:hybrid}
\end{equation}
with ties assigned to the stationary AED. An event is served by a mobile AED if $\tau^{M}_{r,v}(e) < \tau^{S}_{r}(e)$, and receives no response if both terms are infinite, which can occur only at $r = 100\%$, when no stationary AED remains. At $r = 0$, Equation~\eqref{eq:hybrid} reduces to the full stationary network, $T_0(e) = \tau^{S}_{0}(e)$.

\subsection{Performance Measures}
\label{sec:methods-measures}

The improvement of event $e$ is the reduction in response time relative to the full stationary network,
\begin{equation}
\Delta_{r,v}(e) = T_0(e) - T_{r,v}(e).
\label{eq:improvement}
\end{equation}
Each event is classified as improved ($\Delta_{r,v}(e) > 0$), worsened ($\Delta_{r,v}(e) < 0$), unchanged ($\Delta_{r,v}(e) = 0$), or as receiving no response; the four classes partition $\mathcal{E}$ at every ratio. Mean improvement is the average of $\Delta_{r,v}(e)$. %over the events that receive a response. % bahar: in the Toronto P3 runs no event goes unanswered, so it only matters at r = 100 where no-response occurs. 
Reliability at threshold $\theta$ is the share of all events reached within $\theta$ minutes,
\begin{equation}
R_{\theta}(r,v) = \frac{1}{|\mathcal{E}|}\sum_{e\in\mathcal{E}} \mathbf{1}\{T_{r,v}(e) \le \theta\},
\qquad \theta \in \{3, 5, 10\}.
\label{eq:reliability}
\end{equation}
%so that events without a response count as not reached.

%The mobile mode share is the share of events served by a mobile AED. %we don't use the figures related to mode share int he paper.

Because the response rule in Equation~\eqref{eq:hybrid} is a minimum over a random equipped set, every measure is computed from event-level outcomes within each mobility sample and only then averaged across the ten samples. %\textcolor{gray}{averaging distances across samples before the comparison would describe no realizable deployment.} 
Annual values pool all events of 2024, which is equivalent to weighting each month by its number of events, and monthly values are used to show variation across the year. One equipped fleet serves every event within a mobility sample, so outcomes are correlated across events and the effective number of independent draws is the number of samples rather than the number of events. %\textcolor{gray}{; we therefore report the spread across mobility samples and months rather than event-level standard errors.}

\backmatter

% Bahar: Nature Portfolio journals require these two statements. Draft wording; adjust the VFH sentence to the terms of the data-sharing agreement with the City of Toronto.
\bmhead{Data availability}
The Toronto Paramedic Services Incident Data and the City of Toronto AED inventory are available from the City of Toronto Open Data Portal. The New York City EMS Incident Dispatch Data and AED inventory are available from NYC Open Data, and the High-Volume For-Hire Vehicle trip records from the NYC Taxi \& Limousine Commission. Census boundaries and dissemination-area populations are available from Statistics Canada, and road network data from OpenStreetMap. The Toronto Vehicle-for-Hire records were provided by the City of Toronto and are not publicly available.
%; requests for access should be directed to the City of Toronto.

\bmhead{Code availability}
The code used for the analysis is available from the corresponding author upon reasonable request.

\bmhead{Acknowledgments}
The authors thank the City of Toronto Transportation Services Data and Analytics team for providing the Toronto ride-hailing (vehicle-for-hire) mobility data used in this study. This research was granted a research award from the TD Management Data and Analytics Lab (MDAL) at the Rotman School of Management, University of Toronto.

% Appendix here
% Options are (1) APPENDIX (with or without general title) or 
%             (2) APPENDICES (if it has more than one unrelated sections)
% Outcomment the appropriate case if necessary
%
% \begin{APPENDIX}{<Title of the Appendix>}
% \end{APPENDIX}
%
%   or 
%
% \begin{APPENDICES}
% \section{<Title of Section A>}
% \section{<Title of Section B>}
% etc
% \end{APPENDICES}

% The sn-nature document-class option selects the official Nature bibliography style.
\bibliography{references}

\newpage

% =====================================================================
% APPENDIX: Borough-level analysis, NYC (12 months of 2024, FHV fleet)
% Place before \bibliography{references}.
% Requires: booktabs (add \usepackage{booktabs} to the preamble if absent)
% Figures: overlay_F1_mean_improvement.png, overlay_F2_worsened.png,
%          overlay_F3_5min.png, overlay_F4_mode_share.png
% =====================================================================

\begin{appendices}

\section{Data Construction}
\label{app:data}

\subsection{OHCA-Proxy Events}
For Toronto, we use the 2024 Toronto Paramedic Services Incident Data~\citep{TorontoParamedic2025} and retain medical incidents dispatched at the Echo determinant of the Medical Priority Dispatch System, the determinant assigned to suspected cardiac and respiratory arrest, within the city's Forward Sortation Areas (FSAs). From 11 to 31 December 2024 the published records leave the incident type blank for most entries, so Echo incidents in that window are retained at random at a rate of $0.803$, which matches their number to that expected from the same calendar window in 2023 and 2025. The resulting pool contains $8{,}135$ incidents. Because Echo dispatches include calls that do not prove to be arrests, we draw a simple random sample of $4{,}817$ events with a fixed seed. The sample size is anchored to the $4{,}776$ incidents with a cardiac arrest reported by Toronto Paramedic Services for 2025~\citep{TPSAnnual2025}, scaled to 2024 by the ratio of Echo medical incidents in the two years ($8{,}135$ in 2024 and $8{,}065$ in 2025). Each event is placed at the population-weighted centroid of its FSA,
\begin{equation}
\bar{x}_f \;=\; \frac{\sum_{j\in f} p_j\,x_j}{\sum_{j\in f} p_j},
\label{eq:popcentroid}
\end{equation}
where $j$ indexes the 2021 Census dissemination areas within FSA $f$, $p_j$ is the population of dissemination area $j$, and $x_j$ is its representative point. This gives $96$ distinct event locations.

For New York City, we use the 2024 records of the EMS Incident Dispatch Data~\citep{FDNYEMSDispatch} and retain incidents with a cardiac-arrest or resuscitation final call type (ARREST, ARRESC, ARREFC, ARREFT, CPR, CPRC) and an incident disposition code of 83, 92, or 94, which yields $10{,}082$ events. Each event is placed at the centroid of its ZIP code, giving $181$ distinct event locations, and is assigned to its borough as recorded in the dispatch data.

\subsection{Stationary AEDs}
For Toronto, the stationary layer is the City of Toronto public AED inventory~\citep{TorontoAED2025}, which lists one record per device by address. We geocode each address and place devices that an organization registers under a single head-office address, such as transit stations and public beaches, at their individual sites. Of the $1{,}428$ devices listed, we exclude devices located outside the city and $92$ devices that are mounted in or carried by municipal vehicles and crews. The resulting inventory contains $1{,}291$ AEDs at $902$ distinct locations, of which $389$ ($30.1\%$) are additional devices at locations that already host one. For New York City, the stationary layer is the NYC Open Data public AED inventory~\citep{NYCAED2026}, which lists each site with its number of devices; it contains $11{,}138$ AEDs at $7{,}373$ distinct locations, of which $3{,}765$ ($33.8\%$) are additional devices at locations that already host one. In both cities a location is a distinct geocoded point.

\subsection{Ride-Hailing Mobility}
For Toronto, the mobility layer consists of twelve months of Vehicle-for-Hire (VFH) records for 2024 from the two major ride-hailing platforms operating in the city, obtained from the City of Toronto. The records carry a persistent anonymized vehicle identifier and follow each vehicle through the three activity phases. They contain approximately $146{,}600$ distinct vehicles, about $69{,}500$ of which are active in an average month, and approximately $117$ million idle segments, $81$ million en-route segments, and $75$ million occupied trips. Idle segments and occupied trips are routed on the driving network from the location data provided, and an en-route vehicle is held at the point where it accepted the trip. Where a trip endpoint is reported only at the FSA level, it is taken from the adjacent segment of the same vehicle where one exists and is otherwise drawn from observed drop-off points in the reported FSA; a trip with no location at one end is held at its known end. Further detail on the vehicle location data is withheld for privacy reasons.

For New York City, the mobility layer consists of the High-Volume For-Hire Vehicle trip records published by the NYC Taxi \& Limousine Commission (TLC)~\citep{NYCTLCTripData}, which cover every trip served by the major ride-hailing platforms, with $15$ to $20$ million trips per month. Each record reports the pickup and drop-off times, the pickup and drop-off taxi zones ($263$ zones), and the metered trip distance. For each ordered pair of zones, we route seven candidate paths on the driving network between endpoint pairs drawn uniformly within the two zone polygons, and assign each trip the candidate whose length is closest to its metered distance.

\section{Borough-Level Analysis in New York City}
\label{app:nyc-borough}

\begin{figure}[htbp]
    \centering
    \includegraphics[width=0.6\textwidth]{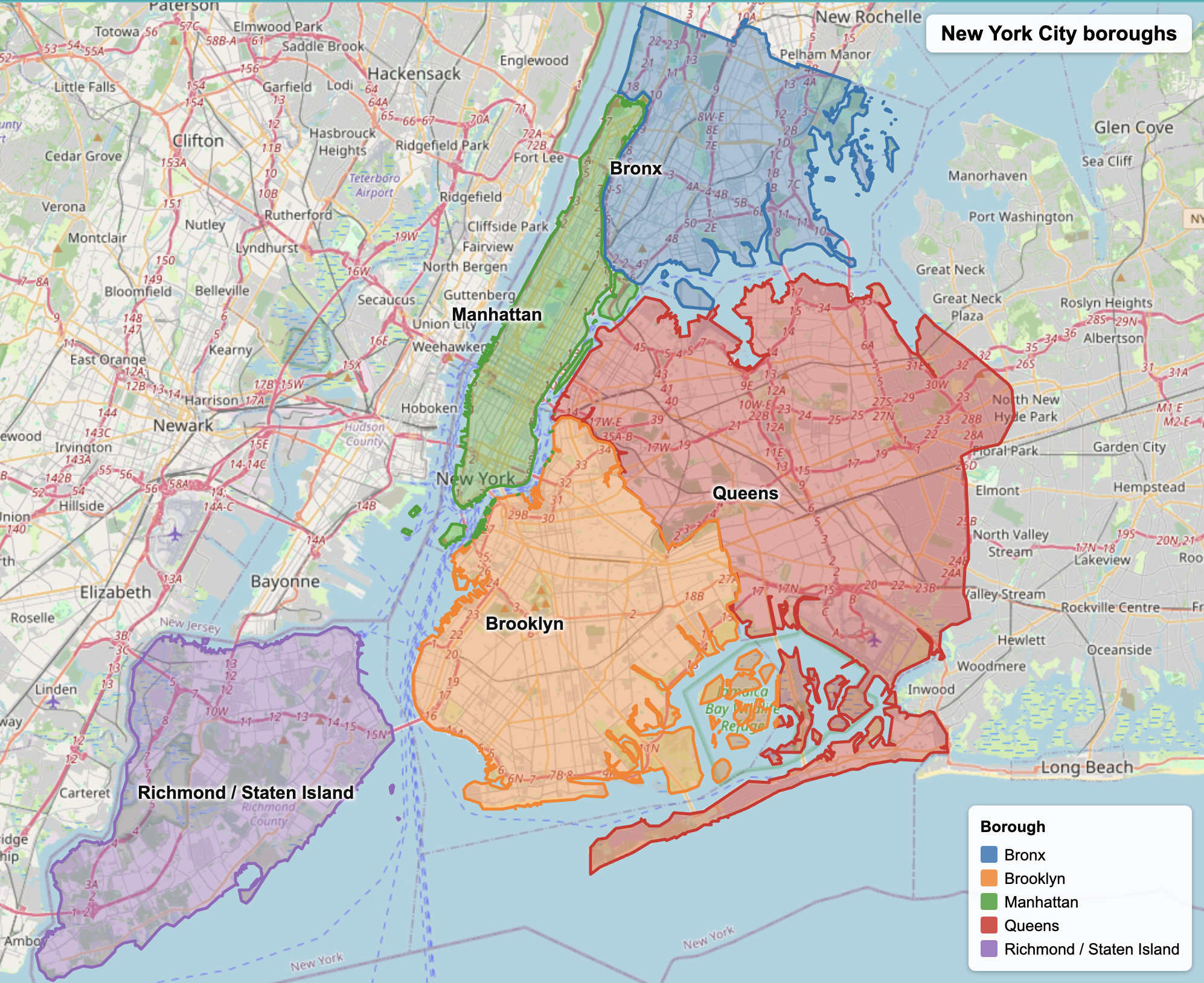}
    \caption{The five boroughs of New York City, shown over a street-level base map. Each borough is shaded in the color used to denote it throughout the borough-level analysis; the boundaries are the official borough boundaries used to assign each OHCA event to the borough in which it occurred.}
    \label{fig:nyc-boroughs}
\end{figure}

The citywide results reported in Section~\ref{sec:toronto-hybrid} for New York City aggregate over an internally heterogeneous urban system, which is itself divided into five boroughs: Manhattan, Brooklyn, Queens, the Bronx, and Staten Island. 
%These boroughs occupy geographically distinct parts of the city %Manhattan is the central island, the Bronx lies directly to its north on the mainland, Brooklyn and Queens together cover the western end of Long Island to the southeast, and Staten Island sits to the southwest across New York Bay, physically separated from the rest of the city. 
This appendix disaggregates the twelve-month 2024 analysis by borough in order to examine how the hybrid policy behaves across urban settings that differ in stationary AED coverage, OHCA geography, and mobility-fleet density.
Figure~\ref{fig:nyc-boroughs} shows the spatial extent of the five boroughs and the boundaries used to assign each OHCA event to a borough; the color assigned to each borough there is carried through the borough-level figures reported below.

Each OHCA event is assigned to the borough in which it occurred. Neither the stationary inventory nor the mobility fleet is restricted to that borough: an event is served by the faster of the nearest remaining stationary AED and the nearest equipped vehicle anywhere in the city, exactly as in the citywide analysis. The borough is therefore a stratification of demand rather than a partition of supply, and the construction of Section~\ref{sec:policy-eval} is otherwise unchanged. Outcomes are classified within each of the ten mobility samples and then averaged across samples, so that each sample represents one realized deployment. The $10{,}082$ OHCA-proxy events of 2024 are distributed unevenly across boroughs, ranging from $2{,}900$ in Brooklyn to $607$ in Richmond/Staten Island.

% \subsection{Baseline Stationary Coverage}
% \label{app:borough-baseline}

Table~\ref{tab:borough_baseline} reports the status quo stationary policy coverage benchmark by borough. Status quo performance differs substantially across boroughs. Manhattan has both the shortest nearest-AED distances and the highest reliability, with $69.4\%$ of OHCA events reachable within $5$~minutes under two-way on-foot retrieval. Richmond/Staten Island has the longest mean distance and the weakest reliability at the $10$-minute threshold, at $36.4\%$.
The Bronx, Brooklyn, and Queens occupy an intermediate range on distance but differ in the shape of their distributions: the mean nearest-AED distance in Queens exceeds its median by over $100$~m, indicating a long right tail of poorly covered events.

\begin{table}[!b]
\centering
\caption{Baseline stationary access by borough under the full stationary network (hybrid ratio $0\%$), New York City, 2024. Distances are network distances to the nearest
AED; times are two-way on-foot retrieval times.}
\label{tab:borough_baseline}
\begin{tabular}{lrrrrrr}
\toprule
Borough & OHCA & Mean & Median & $\leq 3$ & $\leq 5$ & $\leq 10$ \\
        & events & dist.\ (m) & dist.\ (m) & min (\%) & min (\%) & min (\%) \\
\midrule
Bronx                    & 1{,}900  & 333.8 & 294.1 & 3.2  & 16.7 & 61.4 \\
Brooklyn                 & 2{,}900  & 313.0 & 298.9 & 5.0  & 25.1 & 68.4 \\
Manhattan                & 2{,}036  & 188.7 & 171.0 & 30.9 & 69.4 & 93.2 \\
Queens                   & 2{,}639  & 382.9 & 281.1 & 4.4  & 23.9 & 62.8 \\
Richmond / Staten Island & 607      & 410.4 & 461.3 & 18.6 & 28.0 & 36.4 \\
\midrule
All New York City        & 10{,}082 & 316.0 & 271.4 & 10.6 & 32.3 & 68.7 \\
\bottomrule
\end{tabular}
\end{table}

% \subsection{Mean Improvement}
% \label{app:borough-improvement}

\begin{figure}[htbp]
    \centering
    \includegraphics[width=0.7\linewidth]{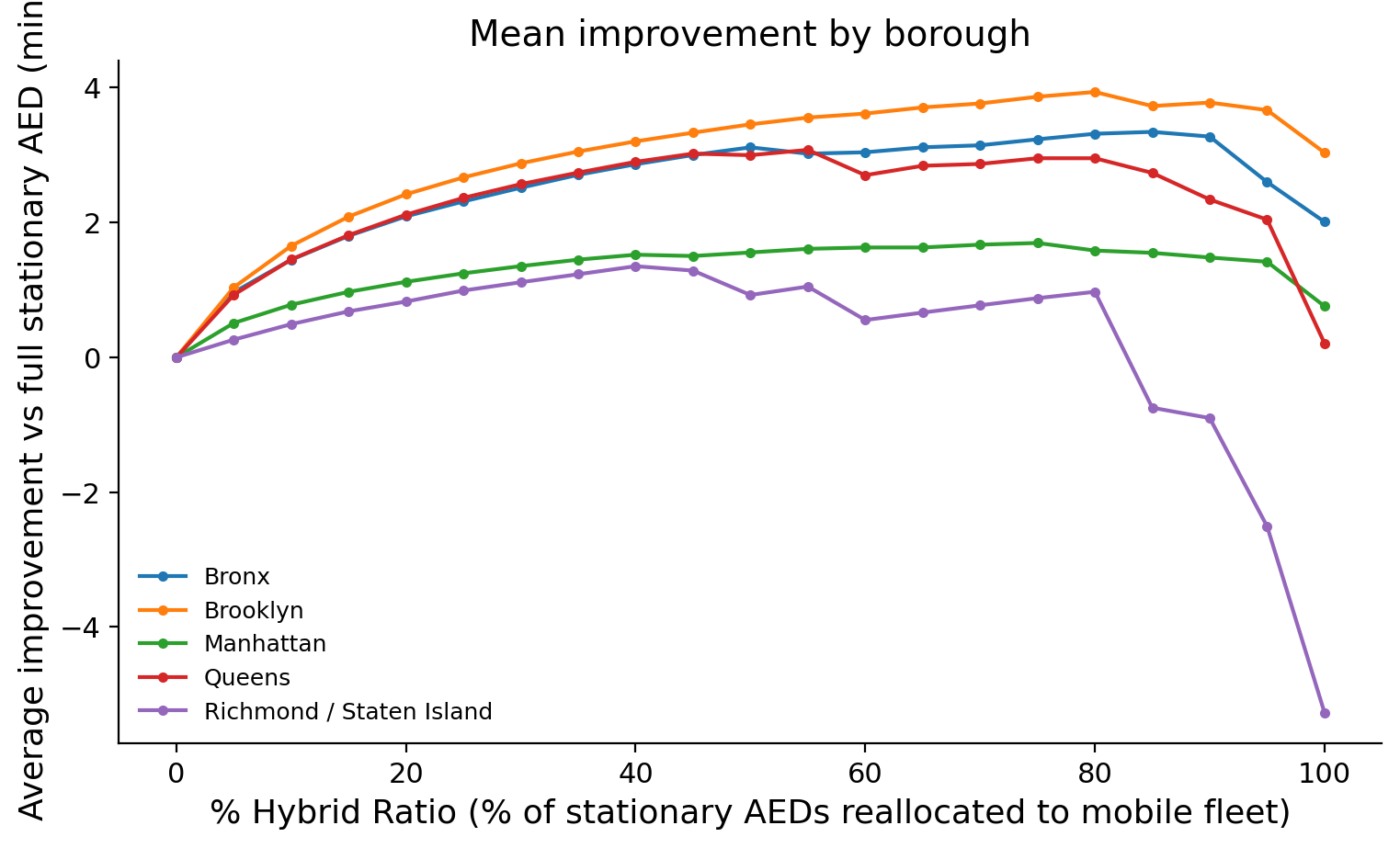}
    \caption{Mean response-time improvement relative to the full stationary network, by borough, as a function of the \emph{hybrid ratio} (New York City, 2024).}
    \label{fig:borough_improvement}
\end{figure}

Figure~\ref{fig:borough_improvement} shows average response-time improvement against the \emph{hybrid ratio} for each borough. All five curves rise steeply at low hybrid ratios, confirming that the early gains observed at the city level are observed in (and driven by) all the boroughs. The curves then diverge in both the location of their maximum and its magnitude.
The interior optimum occurs at a hybrid ratio of $40\%$ in Richmond/Staten Island, $55\%$ in Queens, $75\%$ in Manhattan, $80\%$ in Brooklyn, and $85\%$ in the Bronx, against $80\%$ for the city as a whole. Peak mean improvement ranges from $1.4$~minutes in Richmond/Staten Island and $1.7$~minutes in Manhattan to $3.9$~minutes in Brooklyn.

The ordering of improvement magnitudes is informative: Manhattan, which has the strongest stationary baseline, records the smallest improvement among the four dense boroughs which shows that where a fixed AED is already close, the baseline retrieval time is short and there is little room for the mobile layer to improve on it. Richmond/Staten Island, which has the weakest stationary baseline, also records small improvement, but for the opposite reason, since the mobility fleet is thin there and a dispatched vehicle is frequently distant. The largest gains accrue to boroughs that combine moderate stationary coverage with substantial fleet activity. Improvement depends jointly on the stationary coverage that reallocation removes and on the fleet density that determines what the mobile layer supplies in its place, and interestingly neither factor alone orders the boroughs.

Behavior at high hybrid ratios separates the boroughs further. Under a fully mobile policy, mean improvement remains positive in the Bronx ($2.0$~minutes), Brooklyn ($3.0$~minutes), and Manhattan ($0.8$~minutes), approaches zero in Queens ($0.2$~minutes), and becomes negative in Richmond/Staten Island ($-5.3$~minutes). In that borough, complete reallocation to mobile policy degrades access relative to the status quo stationary policy.

\subsection{Heterogeneous Effect of AED-on-wheels Across NYC Boroughs}
\label{app:borough-breadth}

\begin{figure}[htbp]
\centering
\begin{minipage}[t]{0.49\textwidth}
    \centering
    \includegraphics[width=\textwidth]{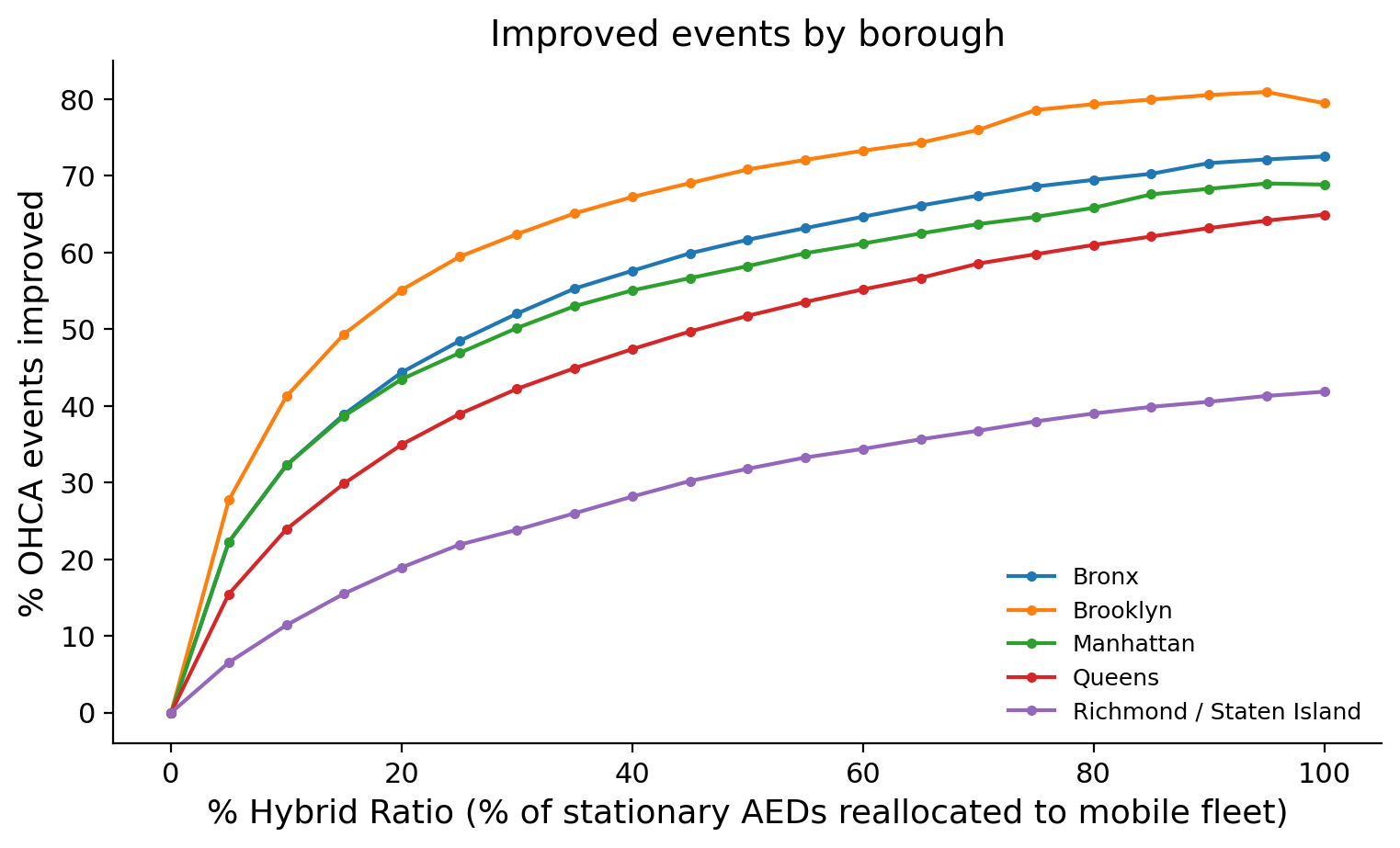}
    \par\smallskip{\small (a) Improved events}
\end{minipage}\hfill
\begin{minipage}[t]{0.49\textwidth}
    \centering
    \includegraphics[width=\textwidth]{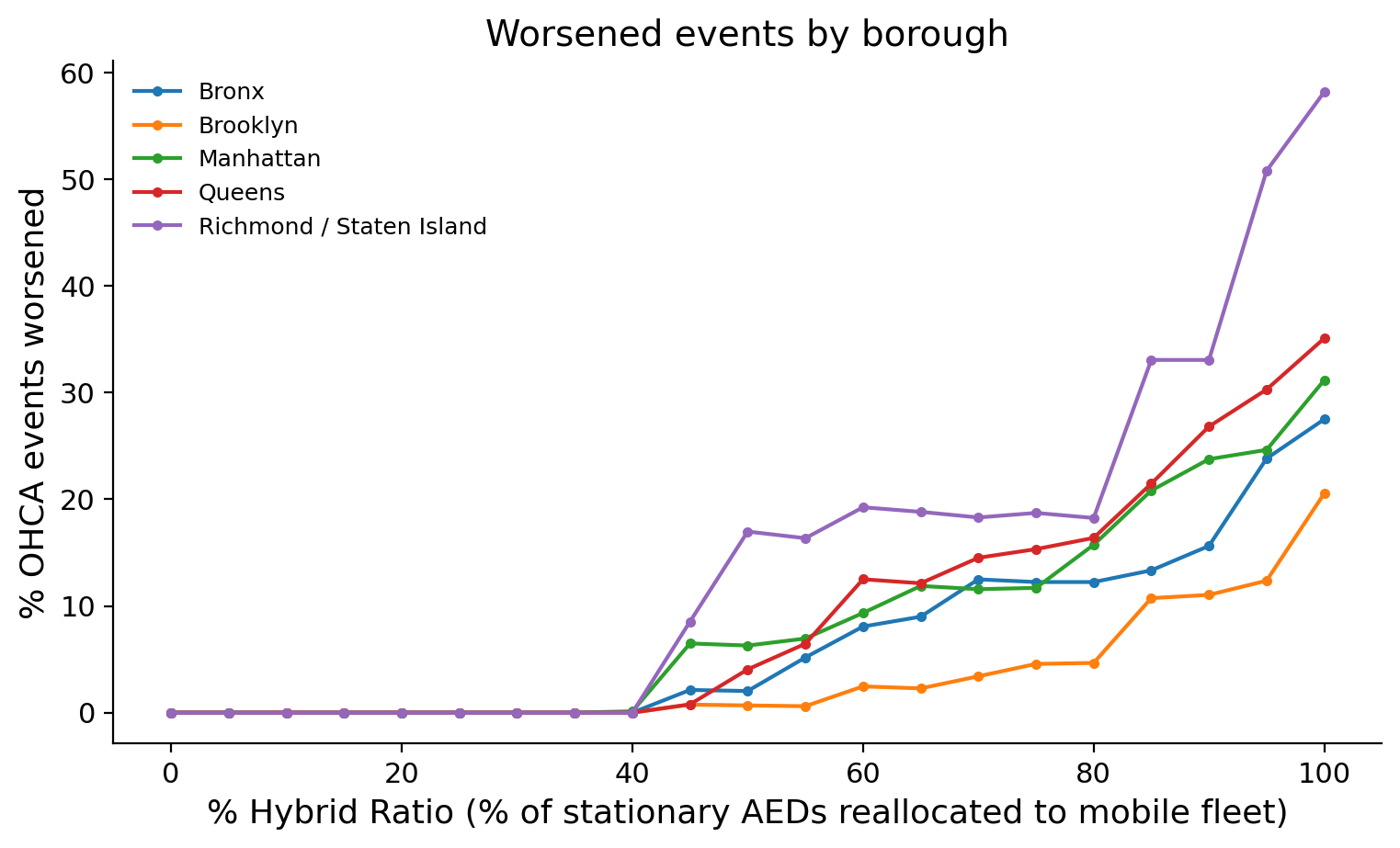}
    \par\smallskip{\small (b) Worsened events}
\end{minipage}
\caption{Share of OHCA events with (a) improved and (b) worsened response time relative to
the full stationary network, by borough, as a function of the \emph{hybrid ratio}
(New York City, 2024).}
\label{fig:borough_outcomes}
\end{figure}

Figure~\ref{fig:borough_outcomes} reports the shares of OHCA events whose response time improves and worsens under the hybrid policy. The improved share rises steeply at low hybrid ratios in every borough and saturates in the upper hybrid ratios without declining: at a hybrid ratio of $20\%$ it spans $18.9\%$ improvement share in Richmond/Staten Island to $55.1\%$ improvement share in Brooklyn, and under a fully mobile policy it spans $41.8\%$ to $79.5\%$. The ordering of boroughs is stable across the entire range, and the vertical separation between the curves
is largest at intermediate hybrid ratios, where boroughs differ most in how quickly mobile AEDs are selected over the stationary alternative.

The worsened share is zero below a hybrid ratio of $40$--$45\%$ in every borough.
This threshold is a consequence of the counterfactual construction rather than an empirical property of mobile AEDs as the stationary rollback removes redundant co-located devices before it removes any site completely, so no fixed AED location disappears until the redundant share of the inventory is exhausted, resulting in no event's stationary access time change. 
Once the threshold is crossed, the worsened share increases monotonically in the hybrid ratio in every borough, but at rates that differ by an order of magnitude. At a hybrid ratio of $60\%$, the worsened share is $2.5\%$ in Brooklyn and $19.2\%$ in Richmond/Staten Island; under a fully mobile policy the corresponding figures are $20.5\%$ and $58.2\%$.

These two panels together explain why mean improvement declines after a certain hybrid ratio while the improved share does not. Beyond the interior optimum, additional reallocation continues to improve access for a slowly growing set of OHCA events, but the events it worsens are losing access by increasingly large margins, and the mean is driven by that magnitude rather than by the count. The interior optimum in each borough is located where the marginal gain from additional mobile supply falls below the marginal loss from withdrawn stationary coverage, and the boroughs in which worsened events accumulate fastest are those whose optima occur at the lowest hybrid ratios.
Table~\ref{tab:borough_summary} summarizes these quantities.

\begin{table}[htbp]
\centering
\caption{Hybrid policy performance by borough, New York City, 2024. The optimum is the hybrid ratio maximizing mean improvement; improved and worsened shares are evaluated at that ratio.}
\label{tab:borough_summary}
\begin{tabular}{lrrrrr}
\toprule
Borough & Optimal & Mean impr. & Improved & Worsened & Mean impr. \\
        & ratio (\%) & at opt.\ (min) & (\%) & (\%) & at $100\%$ (min) \\
\midrule
Bronx                    & 85 & 3.34 & 70.2 & 13.3 & $2.01$ \\
Brooklyn                 & 80 & 3.94 & 79.3 & 4.7  & $3.04$ \\
Manhattan                & 75 & 1.70 & 64.6 & 11.7 & $0.76$ \\
Queens                   & 55 & 3.08 & 53.5 & 6.5  & $0.21$ \\
Richmond / Staten Island & 40 & 1.35 & 28.2 & 0.0  & $-5.27$ \\
\midrule
All New York City        & 80 & 2.91 & 67.5 & 12.2 & $1.14$ \\
\bottomrule
\end{tabular}
\end{table}

\subsection{Reliability}
\label{app:borough-reliability}

\begin{figure}[htbp]
    \centering
    \includegraphics[width=0.7\linewidth]{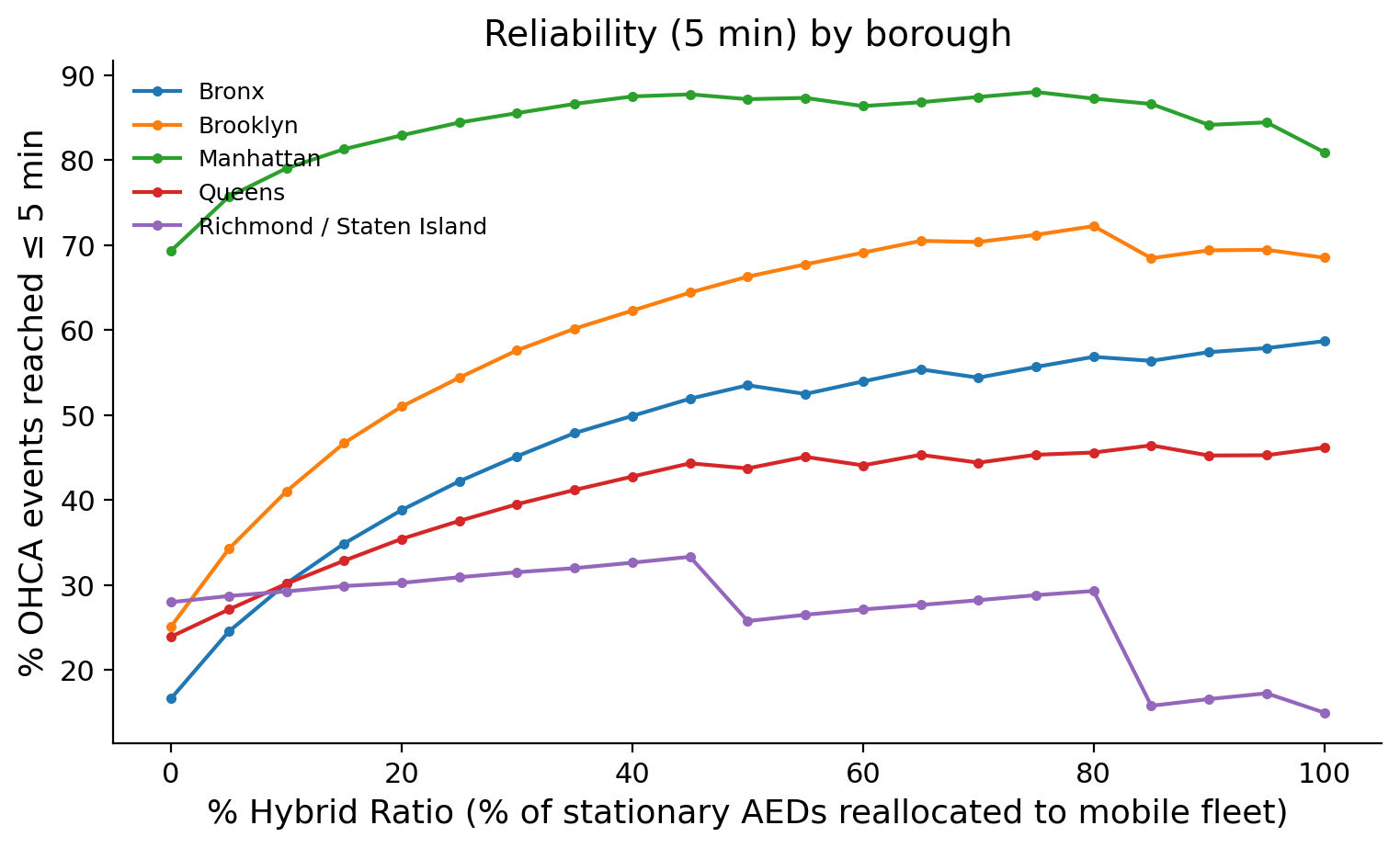}
    \caption{Reliability, measured as the share of OHCA events reached within $5$~minutes,
    by borough, as a function of the \emph{hybrid ratio} (New York City, 2024).}
    \label{fig:borough_reliability}
\end{figure}

Figure~\ref{fig:borough_reliability} reports the share of OHCA events reached within $5$~minutes. Reliability improves in every borough as mobile AEDs are introduced, but the gains are unequal in both size and persistence. Moving from full stationary policy (hybrid ratio 0\%) to fully mobile policy (hybrid ratio 100\%), Brooklyn reliability improves from $25.1\%$ to a maximum of $72.3\%$ and the Bronx from $16.7\%$ to $56.4\%$, while Manhattan improves from an already high $69.4\%$ to $88.0\%$. In Richmond/Staten Island the gain peaks at $32.7\%$ before falling to $15.0\%$ under a fully mobile policy, which is below the stationary baseline. At the $10$-minute threshold, every borough attains its maximum at an interior hybrid ratio between $55\%$ and $85\%$ and declines thereafter, indicating that reliability begins to decline at lower hybrid ratios than mean improvement.
Moreover, the fully mobile policy is dominated on this metric in all five boroughs.
%should I add the 10min figure as well?

% \subsection{Interpretation}
% \label{app:borough-mode}
%and mode share which is commented out!

% Figure~\ref{fig:borough_mode} reports the share of OHCA events served by a mobile AED rather than by a stationary one, which serves as a diagnostic of how frequently the mobility layer is the binding constraint. 
% At a hybrid ratio of $40\%$, the mobile share ranges from $28.2\%$
% in Richmond/Staten Island to $67.2\%$ in Brooklyn. Because the response rule selects the faster of the two options, this share measures how often an equipped vehicle is close enough to displace the stationary alternative, and it therefore reflects local fleet density directly.

% \begin{figure}[htbp]
%     \centering
%     \includegraphics[width=0.8\linewidth]{overlay_F4_mode_share.png}
%     \caption{Share of OHCA events served by a mobile AED, by borough, as a function of the
%     \emph{hybrid ratio} (New York City, 2024).}
%     \label{fig:borough_mode}
% \end{figure}

Taken together, the borough results reproduce, within a single city, the heterogeneity observed also between New York City and Toronto. 
The qualitative structure of the policy is shared: the mean
improvement rises steeply at low hybrid ratios, reaches an interior optimum, and declines; worsened events are absent until the redundant stationary capacity is exhausted and accumulate monotonically thereafter; and reliability at the $10$-minute threshold is maximized at an interior ratio in every case. 
What differs across boroughs is where the optimum falls, how large the attainable gain is, and how sharply performance deteriorates beyond the optimum. %A single citywide hybrid ratio is consequently a coarse instrument: the value that maximizes mean improvement for New York City as a whole, $80\%$, lies well beyond the optimum for Queens and would place Richmond/Staten Island in a regime in which the policy is net-negative on both mean improvement and reliability.

\end{appendices}

\end{document}